\documentclass[12pt,twoside]{article}
\usepackage{times}
\usepackage{amsmath,amssymb}
\usepackage{color}

\allowdisplaybreaks

\def\R {\Bbb R}
\def \N {\Bbb N}

\def\p{\partial}

\def\f{\frac}

\def\la{\lambda}

\def\al{\alpha}
\def\t{\tilde}

\def\vp{\varphi}

\def\th{\theta}
\def\g{\gamma}
\def\G{\Gamma}
\def\si{\sigma}
\def\dl{\delta}

\def\o{\omega}

\def\ds{\displaystyle}
 \allowdisplaybreaks

\begin{document}
 \footskip=0pt
 \footnotesep=2pt
\let\oldsection\section
\renewcommand\section{\setcounter{equation}{0}\oldsection}
\renewcommand\thesection{\arabic{section}}
\renewcommand\theequation{\thesection.\arabic{equation}}
\newtheorem{claim}{\noindent Claim}[section]
\newtheorem{theorem}{\noindent Theorem}[section]
\newtheorem{lemma}{\noindent Lemma}[section]
\newtheorem{proposition}{\noindent Proposition}[section]
\newtheorem{definition}{\noindent Definition}[section]
\newtheorem{remark}{\noindent Remark}[section]
\newtheorem{corollary}{\noindent Corollary}[section]
\newtheorem{example}{\noindent Example}[section]

\title{
On the existence of low regularity solutions to semilinear
generalized Tricomi equations  \\ in mixed type domains}

\author{Ruan,
Zhuoping$^{1, *}$,\qquad Witt, Ingo$^{2}$, \qquad Yin
Huicheng$^{1, }$\footnote{ Ruan Zhuoping and Yin Huicheng were
supported by the NSFC (No.~10931007, No.~11025105),
  and by the Priority
Academic Program Development of Jiangsu Higher Education
Institutions. This research was started when Ruan Zhuoping and
Yin Huicheng were visiting the Mathematical
Institute of the University of G\"ottingen in February-March of 2013.}\vspace{0.5cm}\\
\small 1. Department of Mathematics and
IMS, Nanjing University, Nanjing 210093, China.\\
\vspace{0.5cm}
\small 2.
Mathematical Institute, University of G\"{o}ttingen,
Bunsenstr.~3-5, D-37073 G\"{o}ttingen, Germany. }

\date{}
\maketitle

\centerline {\bf Abstract} \vskip 0.3 true cm
In [19-20], we have established the existence and singularity structures of low regularity solutions to the semilinear
generalized Tricomi equations in the degenerate hyperbolic regions and to the higher order degenerate hyperbolic equations, respectively.
In the present paper,
we shall be concerned with  the  low regularity  solution problem
for the semilinear mixed type equation $\p_t^2u-t^{2l-1}\Delta u= f(t,x,u)$ with an
initial data $u(0,x)=\varphi(x)\in  H^{s}(\Bbb R^n)$ ($0\le s<\f{n}{2}$),
where $(t,x)\in\Bbb R
\times\Bbb R^n$,  $n\ge 2$, $l\in\Bbb N$,
$f(t,x,u)$ is $C^1$ smooth  in its arguments and has  compact support  with respect to the variable $x$.
Under the assumption of the subcritical growth of $f(t,x,u)$ on
$u$, we will show the  existence and regularity of the considered
solution in the mixed type domain $[-T_0, T_0] \times \R^n $ for some fixed constant $T_0>0$.
\vskip 0.3 true cm

{\bf Keywords:} Generalized Tricomi equation, mixed type equation,
confluent hypergeometric function, modified Bessel functions, Calder\'on-Zygmund decomposition,
multiplier  \vskip 0.3 true cm

{\bf Mathematical Subject Classification 2000:} 35L70, 35L65,
35L67, 76N15

\section{Introduction}
In [19-20], we have established the existence and singularity structures of low regularity solutions to the semilinear
generalized Tricomi equations in the hyperbolic regions and to the higher order degenerate hyperbolic equations, respectively.
In the present  paper, we have a further study on the existence and regularities of solutions to the following
$n$-dimensional semilinear generalized Tricomi equation in the mixed type domain
$\Bbb R\times\Bbb R^n$
\begin{equation}
\left\{
\begin{aligned}
&\p_t^2u-t^{2l-1}\Delta u= f(t,x,u),\qquad (t,x)\in \Bbb R\times\Bbb R^n,\\
&u(0,x)=\varphi(x),\qquad\qquad \qquad \qquad  x\in\Bbb R^n,
\end{aligned}
\right.\tag{1.1}
\end{equation}
where $l\in\Bbb N$, $x=(x_1, ..., x_n)$, $n\ge 2$,
$\Delta=\ds\sum_{i=1}^n\p_i^2$, $\varphi(x)\in  H^s(\Bbb R^n)$ ($0\le s<\ds\f{n}{2}$),
$f(t,x,u)$ is $C^1$ smooth in
its arguments and has a compact support $E$ on the variable $x$. Moreover,
for any $T>0$, there exists $C_T >0$ such that
for $(t,x,u)\in [-T, T]\times  E\times\Bbb R$,
$$|f(t,x,u)|\le C_T(1+|u|)^{\mu}\quad\text{and\quad $|\p_uf(t,x,u)|\le C_T(1+|u|)^{\max\{\mu-1, 0\}}$ },\eqno{(1.2)}$$
where $C_T>0$ is a constant depending only on $T$,  and the exponent $\mu\ge 0$ fulfills
$$ \mu < p_0\equiv\f{2n}{n-2s}.\eqno{(1.3)} $$
Here we point out that the number $p_0$ defined in (1.3) comes from the
Sobolev imbedding formula $H^{s}(\Bbb R^n)\subset L^{p_0}(\Bbb R^n)$. Thus (1.2) and (1.3) mean that
the nonlinearity $f$ in (1.1) admits a ``subcritical'' growth on the variable $u$. In addition,
we shall illustrate that the scope of the exponent $\mu$ for
solving the problem (1.1) is closely related to the number $Q_0\equiv 1+\ds\f{n(2l+1)}{2}$.
In the terminology of [16] and the references therein, $Q_0$ is called the homogeneous dimension corresponding to the degenerate elliptic operator $\p_t^2u-t^{2l-1}\Delta$
for $t\le 0$. Our  main result in this paper is:

\vskip .2cm

{\bf Theorem 1.1.} {\it Under the conditions (1.2)-(1.3), there exists a constant $T_0>0$ such that the problem (1.1)
has a unique solution $u(t,x)\in C([-T_0, T_0], L^{p_0}(\Bbb R^n))$ when $0 \le \mu \le 1$, or when $1 < \mu <p_0$ and
$Q_0 \le \ds\f {p_0}{\mu-1}$.}

\vskip .2cm

{\bf Remark 1.1.} {\it
It seems to be necessary that the restriction
of $Q_0 \le \ds\f {p_0}{\mu-1}$ is posed when $1 < \mu <p_0$ in Theorem 1.1, otherwise,
the standard iteration scheme for solving the problem (1.1) only works in finite steps
or the solution $u\not\in C([-T_0, T_0], L^{p_0}(\Bbb R^n))$ when
$Q_0>\ds\f {p_0}{\mu-1}$. One can see Remark 4.1
and the related explanations below (4.7) in $\S 4$, respectively.}

\vskip .2cm

{\bf Remark 1.2.} {\it For $l=1$, (1.1) is the well-known semilinear Tricomi equation $\p_t^2u-t\Delta u=
f(t,x,u)$. When an initial data $u(0,x)=\varphi(x)\in H^s(\Bbb R^n)$  with $s>\ds\f{n}{2}$
is given and the {\bf crucial assumption} of $supp f\subset\{t\ge 0\}$ is posed (namely,
$f\equiv 0$  holds in $t\le 0$, which means that the related Tricomi equation is {\bf linear} in the elliptic region $\{t\le 0\}$),
M. Beals
in [2] shows that the problem (1.1) has a regular solution $u\in C((-\infty, T], H^s(\Bbb R^n)\cap C^1((-\infty, T],
H^{s-\f{5}{6}}(\Bbb R^n))\cap C^2((-\infty, T],
H^{s-\f{11}{6}}(\Bbb R^n))$  for some constant $T>0$. Here we point out that  the key assumption of $supp f\subset\{t\ge 0\}$
in [2] has been removed as well as the local existence of low regularity solutions is
established in our present paper.}

\vskip .2cm
{\bf Remark 1.3.} {\it If the initial data $\varphi(x)\in L^{\infty}(\Bbb R^n)\cap H^s(\Bbb R^n)$
for $n=2, 3$ and $s\ge 0$ is given, then we can remove all the assumptions in (1.2)-(1.3). In fact, in this case,
from the proof procedure in $\S 4$,
we can derive that the solution $u(t,x)$ of (1.1) satisfies:
$u(t,x),\p_tu(t,x)\in L^{\infty}([-T_0, 0]\times\Bbb R^n)$ for some fixed $T_0>0$. Based on this, together with
the $L_{loc}^{\infty}([0, +\infty)\times\Bbb R^n)$ property of the solution $v(t,x)$ to the
linear equation $\p_t^2v-t^{2l-1}\Delta v=0$ with $(v(0,x), \p_tv(0,x))=(\vp_0(x), \vp_1(x))$, where
$\vp_0(x), \vp_1(x)\in L^{\infty}(\Bbb R^n)$, we may show that (1.1) is also solvable in the degenerate hyperbolic region
in the space $L^{\infty}([0, T_1]\times\Bbb R^n)$ for some positive constant $T_1>0$. Combining these
two cases, we know that (1.1) is locally solvable (see more detailed explanations in Remark 5.1 of $\S 5$ below}).

\vskip .2cm

{\bf Remark 1.4.} {\it If the nonlinearity $f(t,x,u)$ in (1.1) increases sublinearly or linearly with respect to
the variable $u$, namely, $0\le\mu\le 1$ holds
in (1.2), then by a minor modification on the proof of Theorem 1.1, we can obtain a global solution
$u\in C(\Bbb R, L^{p_0}(\Bbb R^n))$ to (1.1). }

\vskip .2cm

{\bf Remark 1.5.} {\it In order to guarantee the existence and uniqueness of the solution $u$
to (1.1), the assumption that $f(t,x,u)$ is compactly supported on
the variable $x$ is required.
Otherwise, in case $f(t,x,u)\sim |u|^{p-1}u$ with $p\ge 1$, due to the the existence
of infinite eigenvalues for the Tricomi operator $\p_t^2
-t\Delta$ with $t\in\Bbb R$, (1.1)  can admit infinitely many solutions or have no solution (one can see
the references [13], [15-16] and so on).}

\vskip .2cm

{\bf Remark 1.6.} {\it With respect to the  existence,  singularity structures,
and singularity propagation theories of classical
solutions to the semilinear generalized Tricomi
equations in the {\bf degenerate hyperbolic regions/mixed type regions} or to the higher order {\bf degenerate hyperbolic} equations,
so far there have been many interesting results (one can see [1-5], [7], [19-20], [25-26] and the references therein).
For the linear Tricomi equation in the mixed type region, when
the closed boundary value is given, the authors in [14] establish the existence and uniqueness of weak solutions.
Here our focus in Theorem 1.1 is on the existence of the low regularity solution to the semilinear problem
(1.1) in the mixed type region $\Bbb R\times\Bbb R^n$.}

\vskip .2cm

{\bf Remark 1.7.} {\it For the linear Tricomi equation $\p_x^2u-x\p_y^2u=0$ with an initial data
$u(0,y)=u_0(y)$, its solvability in the whole region $\Bbb R^2$ has some very important practical applications in
the continuous transonic gas dynamics of isentropic and irrotational flows. Indeed, once the solution $u$ is found,
one can seek out the corresponding de Laval nozzle walls through its streamlines and then the required
position of the sonic curve (corresponding to the lone $\{x=0\}$ in the hodograph plane) in the
de Laval nozzle can be determined (more introductions on this physical background can be found in [12] and so on).}

\vskip .2cm

For $n=1, l=1$ and $f(t,x,u)\equiv 0$, the equation
in (1.1) becomes the classical
Tricomi equation which arises in transonic gas dynamics and has been
extensively investigated in bounded domain with suitable boundary
conditions from various viewpoints (one can see the review paper [18] and the references therein).
For $l=1$ and $n=2$, with respect to
the equation $\p_t^2u-t\triangle u=f(t,x,u)$ with an
initial data $u(0,x)\in H^s(\Bbb R^n)$ $(s>\ds\f{n}{2})$,
under the crucial assumption of $supp f\subset\{t\ge 0\}$, M. Beals in [2] established the local
existence of the solution $u(t,x)\in C((-\infty, T], H^s(\Bbb R^n))\cap C^1((-\infty, T],
H^{s-\f{5}{6}}(\Bbb R^n))\cap C^2((-\infty, T], H^{s-\f{11}{6}}(\Bbb
R^n))$ for some $T>0$, moreover the  $H^s(\Bbb R^n)$ conormal regularity of
$u$ with respect to the characteristic surfaces $x_1=\pm
\ds\f{2}{3}t^{\f{3}{2}}$ was also obtained in [2]. For more
general nonlinear degenerate hyperbolic equations with discontinuous initial data, the authors in [19-20] obtained the local
existence  of low regularity solutions. In the
present paper, we focus  on the {\bf low regularity solution} problem for
the semilinear generalized Tricomi equation with an initial data $u(0,x)$ in the
{\bf mixed type region $\Bbb R\times\Bbb R^n$}.

We now comment on the proof of Theorem 1.1. In order to establish the
existence and regularity of the solution to (1.1) with an initial data
$u(0,x)=\varphi(x)$, we first consider the linear equation $\p_t^2v-t^{2l-1}\Delta v=0$ with
$v(0,x)=\varphi(x)$  in the domain $\Bbb R\times\Bbb R^n$
and obtain $v(t,x)\in C(\Bbb R, H^{s})\cap C^1(\Bbb R, H^{s-\f{4l+1}{2(2l+1)}})$. Subsequently, we set $w=u-v$ and
get a second order nonlinear degenerate equation of $w$ from (1.1) in the domain $\{t\le 0\}$.
By utilizing some delicate harmonic analysis methods (e.g., Calder\'on-Zygmund decomposition, interpolation,  multiplier, fractional integral,
and so on) as in [8-9], we can establish some suitably weighted
$W^{2,p}(\Bbb R^+\times\Bbb R^n)$ estimates on $w$ and further obtain the solvability of $w$ in $\{t\le 0\}$ by the fixed point principle.
From this, one can get another initial data $\p_t u(0,x)$, which is necessary to  solve
(1.1) in the degenerate hyperbolic region $\{t\ge 0\}$. Finally by using some techniques
in degenerate hyperbolic equations (see [19-20], [25-26] and the references therein),
we  can establish the solvability and regularity of the solution $u$ to (1.1) in the domain $\{t \ge 0\}$.
Then a local solution $u$ in the mixed type domain $\Bbb R\times\Bbb R^n$
could be obtained by patching  the two solutions
got in the degenerate elliptic domain and the degenerate hyperbolic domain separately.

This paper is organized as follows. In $\S 2$, we will give some preliminary results
and useful estimates on the solutions to the linear degenerate elliptic equation
$\p_t^2u+t^m\Delta u=t^mf(t,x)$ with $u(0,x)=0$
and $f(t,x)\in L^p(\Bbb R^+\times\Bbb R^n)$, where $m\in\Bbb N$ and $t\ge 0$.  In
$\S 3$, we establish more general weighted $W^{2,p}$ estimates of the solutions
to the equation  $\p_t^2u+t^m\Delta u=t^\nu f(t,x)$ in  $\{t\ge 0\} $ for $0<\nu<m$.
Here we emphasize that such weighted estimates also admit independent interests in
the degenerate elliptic equations.
Based on the results in $\S 2$ and $\S 3$, we can show the local
existence and regularity of (1.1) in the degenerate elliptic region $\{t\le 0\}$ in $\S 4$.
Moreover, we may obtain $\p_tu(0,x)$ from Theorem 4.1 in $\S 4$. Together with the initial data $u(0,x)$, we
can solve (1.1)
locally in the hyperbolic region $\{t\ge 0\}$ in $\S 5$, and subsequently complete the proof of Theorem 1.1
in $\S 6$.

\vskip 0.4 true cm
\section{Some preliminaries and $W^{2,p}$ estimates for the inhomogeneous generalized Tricomi  equation  }

In this section, we mainly study the $W^{2,p}$ regularity of the solution $u(t,x)$ to the linear generalized Tricomi equation
$\p_t^2u+t^{m}\Delta u=t^m g(t,x)$
with the boundary value $u(0,x)=0$ for  $t\ge 0$, $m\in\Bbb N$ and $g(t,x)\in L^p(\Bbb R^+\times\Bbb R^n)$ $(1<p<\infty)$.
To this end, we require to apply some harmonic analysis tools (e.g., Calder\'on-Zygmund decomposition,
generalized H\"ormander's multiplier theorem and so on) and some properties of  modified Bessel functions.
The so-called modified Bessel function  $K_\nu(t)= \int_0^\infty e^{-t
\cosh z } \cosh (\nu z) dz$ $(\nu\in\Bbb R)$ is a solution to the equation $\Big(t^2\ds \f {d^2} {dt^2} +
t\f {d}{dt} - (t^2 + \nu^2) \Big) K_\nu(t)=0$ for $t>0$, moreover, there holds that $\ds\lim_{t\to+\infty}K_\nu(t)=0$
and $\ds\lim_{t\to 0+}t^{\f12}K_\nu(t)=C_\nu>0$ (one can see more properties of $K_\nu(t)$ in [6], [23] and so on). As in [8-9], set
$\lambda(t)=\ds C_{\f{1}{m+2}} t^{\f 1 2} K_{\f 1 {m+2}}(\f 2 {m+2} t^{\f {m+2} 2})$ for $t>0$ and $m\in\Bbb N$.
Then a direct verification yields

\vskip .3cm

{\bf  Lemma 2.1.}  {\it For $t\ge 0$,

(i) $\la (t)$ is a solution to the equation
$u''(t) -t^m u(t)=0 $ with $u(0)=1$ and $u(+\infty)=0$. \vskip .1cm

(ii) The equation $u''(t)
-t^m u(t)=g(t)$ with $u(0)=1$ and $u(+\infty)=0$ has a solution
$$u(t)= - \la(t) \int_0^\infty \Big( \int_0^{ \min(t,\sigma)} \f {1} {\la^2(y)} dy \Big) \la(\sigma)g(\sigma) d\sigma.$$
}

We now cite the Lemma 2.1 of [8], which illustrates some basic properties of $\la(t)$.
\vskip .3cm

{\bf Lemma 2.2.} {\it

(i) $\la(t)$ and $-\la'(t)$  are decreasing. \hfill (2.1) \vskip .1cm

(ii) $\la(t) + |\la'(t)| \le C_M (1+t)^{-M}$ holds for any $M \in \N$, where $C_M$ is a positive constant
depending on $M$. \hfill (2.2)
 \vskip .1cm

(iii) $\ds\f {\la(a)} {\la (b)} \le ( \f a b )^{\f 1 2} exp \Big( \f 2 {m+2}
b^{\f {m+2} 2} - \f 2 {m+2} a^{ \f {m+2} 2}\Big) $ when $a \ge b
>0.$ \hfill (2.3)
\vskip .1cm

(iv) There exists a constant $C>1$ such that $\ds\f 1 C  \le \f {|\la'(t)|}  {\la(t) t^{m/2}}$ for all $t>0$ and
$\ds\f {|\la'(t)|} {\la(t) t^{m/2}} \le C$ for $t\ge 1$ hold.  \hfill (2.4)}

\vskip .2cm

Next we prove the  global existence and regularity of the solution to the linear generalized Tricomi equation
with an initial data in the whole mixed-type domain $\Bbb R\times\Bbb R^n$.

\vskip .3cm

{\bf Lemma 2.3.} {\it Consider the problem
\begin{equation}
\left\{
\begin{aligned}
&\p_t^2v-t^{2l-1}\Delta v=0,\qquad (t,x)\in \Bbb R\times\Bbb R^n,\\
&v(0,x)=\psi(x),\qquad  \qquad   x\in\Bbb R^n,
\end{aligned}
\right.\tag{2.5}
\end{equation}
where $l\in\Bbb N$, $\psi(x)\in H^{\gamma}(\Bbb R^n)$ $(\gamma\in\Bbb R)$, then (2.5)
has a solution $v(t,x)\in C(\Bbb R, H^{\gamma})\cap C^1(\Bbb R, H^{\gamma-\f{2l+3}{2(2l+1)}})$.
}

\vskip .2cm
{\bf Proof.} At first we study (2.5) in the elliptic region $\{t\le 0\}$
\begin{equation}
\left\{
\begin{aligned}& \p_t^2 w - t^{2l-1}\Delta w=0,  \qquad   (t,x)\in (-\infty,0] \times
\R^n,\\
&w(0,x)=\psi(x).
\end{aligned}
\right.\tag{2.6}
\end{equation}
Taking
Fourier transform with respect to the variable $x$ in (2.6) yields
\begin{equation}
\left\{
\begin{aligned} & \p_t^2 \hat{w}(t,\xi) + t^{2l-1}|\xi|^2 \hat{w}(t,\xi)
=0, \qquad (t, \xi)\in  (-\infty,0] \times
\R^n,\\
& \hat{w}(0,\xi) = \hat{\psi}(\xi).
\end{aligned}
\right.\tag{2.7}
\end{equation}
By Lemma 2.1.(i), the solution
$\hat{w}(t,\xi)$ to (2.7) can be expressed as
$$\hat{w}(t,\xi)=\la(-ts) \hat{\psi}(\xi)\quad \text{with $s= |\xi|^{\f 2 {2l+1}}$}.\eqno{(2.8)}$$
Then it follows from (2.8) and Lemma 2.2.(ii)  that
\begin{align*}\|w(t,\cdot)\|_{H^\gamma(\R^n)}=&  \| \la(-ts)
<\xi>^\gamma
\hat{\psi}(\xi) \|_{L^2(\R^n)}  \\
\le & C \|<\xi>^\gamma
\hat{\psi}(\xi)\|_{L^2(\R^n)} =C
\|\psi\|_{H^\gamma(\R^n)}\tag{2.9}
\end{align*}

and
\begin{align*} \|\p_t w(t,\cdot)\|_{H^{\gamma-\f{2}{2l+1}}(\R^n)}=&  \| \la'(-ts)
s<\xi>^{\gamma-\f{2}{2l+1}}
\hat{\psi}(\xi) \|_{L^2(\R^n)}\\
&\le C \|s<\xi>^{\gamma-\f{2}{2l+1}}
\hat{\psi}(\xi)\|_{L^2(\R^n)} \\ & \le
C \|\psi\|_{H^\gamma(\R^n)}.\tag{2.10}
\end{align*}

Thus, we have from (2.9) and (2.10)
$$w(t,x)\in C((-\infty,0], H^\gamma(\R^n))\cap C^1((-\infty,0], H^{\gamma-\f{2}{2l+1}}(\R^n)).\eqno{(2.11)}$$

Next we consider the corresponding degenerate hyperbolic part of (2.5) in the region $\{t\ge 0\}$
\begin{equation}
\left\{
\begin{aligned} & \p_t^2 u- t^{2l-1}\Delta u=0,  \qquad (t,x)\in
[0, +\infty)\times
\R^n,\\
&u(0,x) = \psi(x), \qquad \p_t u(0,x) = \p_t w(0,x),
\end{aligned}
\right.\tag{2.12}
\end{equation}
where $\p_t w(0,x)\in H^{\gamma -\f {2}
{2l+1} }(\R^n)$ comes from (2.11).

Upon applying Proposition 3.3 in [19], we arrive at
$$u(t,x)\in C([0,
+\infty), H^\gamma(\Bbb R^n ))\cap C^1([0, \infty),
H^{\gamma-\f{2l+3}{2(2l+1)}}(\Bbb R^n)).\eqno{(2.13)}$$
Combining (2.11) with (2.13) yields that  problem (2.5)
has a solution $v(t,x)$ satisfying
$$ v(t,x) \in C(\Bbb R, H^{\gamma}(\R^n))\cap C^1(\Bbb R,
H^{\gamma-\f{2l+3}{2(2l+1)}}(\Bbb R^n)),$$
therefore, we  complete the proof of Lemma 2.3. \hfill $\square$

\vskip .2cm

To get the solvability of  the nonlinear problem (1.1),  as the first step we  intend to solve (1.1) in the degenerate elliptic region
$(-\infty, 0]\times\Bbb R^n$. To do so, we require to derive the weighted $W^{2,p}$ estimate of the solution to the problem   $\p_t^2 w+t^{m}\Delta w=t^\nu g(t,x)$ ($0\le \nu\le m$)
with $w(0,x)=0$ in $\{t\ge 0\}$ so that (1.1) can be solved by applying the
Hardy's inequality and the fixed point theorem (one can  see the details in $\S 4$
below), where $g(t,x)\in L^p(\Bbb R^{n+1}_+)$.   In this section, we only treat  the case of $\nu=m$,
where the $W^{2,p}$ (not weighted $W^{2,p}$) estimate can be derived. Based on this,  by the interpolation method we can
obtain the weighted $W^{2,p}$ estimates for the general $\nu$ (see Theorem 3.1 and its proof in $\S 3$ below).

From the equation $\p_t^2w+t^{m}\Delta w=t^m g(t,x)$
with $w(0,x)=0$, we have
\begin{equation}
\left\{
\begin{aligned} &\p_t^2 \hat{w}(t,\xi)- t^{m}|\xi|^2 \hat{w}(t,\xi)=t^m \hat{g}(t,\xi), \quad (t,\xi) \in [0,+\infty) \times\Bbb R^n,\\
&\hat{w}(0,\xi)=0.
\end{aligned}
\right.\tag{2.14}
\end{equation}
This, together with Lemma 2.1.(ii), yields
$$\hat{w}(t,\xi)=\int_0^{+\infty} \hat{T}(t, \sigma, \xi)  \sigma^m \hat{g}(\sigma,\xi) d\sigma, \eqno{(2.15)}$$
where  $ \hat{T}(t, \sigma, \xi)=\ds\int_0^{\min(t, \sigma)} \f
{\la(ts) \la(\sigma s)} {\la^2(ys)} dy$ and $s= |\xi|^{2/{m+2}}$. By (2.15), we have
$$
\Delta w=  \mathcal F_\xi^{-1} \big(|\xi|^2 \hat{u}(t,\xi)
\big)=\mathcal F_\xi^{-1} \Big( \int_0^{+\infty} \hat{K}(t, \sigma,\xi)
\hat{g}(\sigma,\xi) d\sigma \Big),\eqno{(2.16)}$$
here $\hat{K}(t, \sigma, \xi)=|\xi|^2\sigma^m
\hat{T}(t,\sigma,\xi)$. We start to analyze the property of kernel $\hat{K}(t, \sigma, \xi)$
so that the $L^p(\Bbb R^+\times\Bbb R^n)$ estimate of $\Delta w$ can be obtained. To this end,
we require to apply a basic result on the $L^p$ boundedness for a
class of integral operators, which is established in  Theorem 1.1 of [8] (a generalized H\"ormander multiplier
theorem in Theorem 7.95 in [10]). Suppose that
the temperate distribution $K(t, \sigma, x)\in  S'(\R^+ \times
 \R^n)$ ($t\in\Bbb R^+$ is taken as a parameter) satisfies

(i) For each fixed $(t,\si)\in\Bbb R^+\times\Bbb R^+$, $\hat{K}(t, \sigma, \xi)=\mathcal{F}_x\bigl(K(t,\sigma,\cdot)\bigr)(\xi)
\in C^\infty(\R^n\setminus\{0\}) \cap
C(\R^n)$, and $\hat{K}$ is piecewise continuous in $(t, \sigma)$. $\hfill (2.17)$

(ii) $\ds\sup_{t, \xi} \int_{0}^{+\infty} |\hat{K}(t, \sigma, \xi)| d\sigma \le C.\qquad \qquad \qquad
\qquad \qquad \qquad \qquad \qquad \qquad  \hfill (2.18)$

(iii) $\ds\sup_{ \sigma, \xi} \int_{0}^{+\infty} |\hat{K}(t, \sigma, \xi)| dt \le C. \hfill  (2.19)$

(iv) Denote by $\Delta(a,b)=(a-b, a+b)$ and $\mathcal C \Delta(a, 2^{q+1}b)= \R^+\setminus \Delta(a, 2^{q+1}b)$
for some integer $q$
with $q>n$. For $h(\si) \in C^\infty_0(\R^+)$ and $r>0$, there exists a constant $C_q>0$ depending only
on $q$ such that for all $|\al|\le q$

$$\sup_{\f r 2 \le |\xi|\le 2 r} r^{|\alpha|} \int_{0}^{+\infty} \Big|\partial_\xi^\alpha
\int_{0}^{+\infty} \hat{K}(t, \sigma, \xi) h(\sigma) d\sigma \Big| dt \le C_q \int_0^{+\infty}
|h(\sigma)| d\sigma, \eqno{(2.20)}$$

$$\sup_{\f r 2 \le |\xi|\le 2 r} r^{|\alpha|} \int_{\mathcal C \Delta(a, 2^{q+1}b)} \Big|\partial_\xi^\alpha
\int_{0}^{+\infty} \hat{K}(t, \sigma, \xi) h(\sigma) d\sigma \Big| dt \le C_q
D_1(a,b, r)\int_0^{+\infty} |h(\sigma)| d\sigma,  \eqno{(2.21)}$$

and if $  \int_0^{+\infty} h(\sigma)
d\sigma =0,$ then
$$\sup_{\f r 2 \le |\xi|\le 2 r} r^{|\alpha|} \int_{0}^{+\infty} \Big|\partial_\xi^\alpha
\int_{0}^{+\infty} \hat{K}(t, \sigma, \xi) h(\sigma) d\sigma \Big| dt \notag \le C_q
D_2(a,b, r)\int |h(\sigma)| d\sigma,  \eqno{(2.22)} $$
where $D_i(a,b, r)
(i=1,2)$ are some positive functions of $a, b, r$. Then one has

\vskip .1cm

{\bf Lemma 2.4.} ({See Theorem 1.1 of [8]})  {\it  Let (2.17)-(2.22) be fulfilled. Assume that for
each  compactly supported $f(t,x) \in L^1(\R^+\times\Bbb R^{n})$, and for each $\dl>0$,
there is a Calder\'on-Zygmund decomposition
$$f=f_0 + \ds\sum_{k=1}^{\infty} r_k,$$
where $\|f_0\|_{L^1} + \sum \|r_k\|_{L^1} \le 3 \|f\|_{L^1}$, and $\sup |f_0| \le C \dl$. In addition,
for some disjoint cubes $Q_k=\Delta_k(a_k, \tilde{b}_k) \times
I_k(x_k, b_k)$, here $\Delta_k$ (or $I_k$) is the interval (or cube) centered at $a_k$ (or $x_k$)
with side-length $2 \tilde{b}_k$ (or $2 b_k$), assuming that
$$r_k \text{ is supported in } Q_k, \ \ \int_{Q_k} r_k dtdx=0,  \ \ \text{ and }\quad \dl \sum |Q_k| \le \|f\|_{L^1}.$$
Moreover, for any $k\in\Bbb N$, we also suppose that
$$\sum_{j\ge [-\log_2 b_k]} D_1(a_k, \tilde{b}_k, 2^j) + \sum_{j\le [-\log_2 b_k]} D_2(a_k, \tilde{b}_k, 2^j) \le C, \eqno{(2.23)}$$
where $C>0$ is a generic
positive constant independent of $\dl, f, k, a_k, x_k, b_k$ and $\tilde{b}_k$.

Then the operator $K$ defined by $$(Kf)(t,x)= \int_0^\infty K(t, \sigma, \cdot)
\ast f(\sigma,\cdot) d\sigma$$ is bounded on $L^p(\R^+\times\R^{n})$ for all
$p \in (1, \infty)$.}

Next we apply Lemma 2.4 to establish the $L^p$ boundedness of $\Delta w$ in (2.16). For this purpose,
we require to verify  that the kernel $\hat{K}(t, \sigma, \xi)=|\xi|^2\sigma^m
\hat{T}(t,\sigma,\xi)$ in (2.16) satisfies (2.17)-(2.22). Here we point out that
our analysis for $\hat{K}(t, \sigma, \xi)$ is much more delicate and involved than that for the kernel $\hat{K_0}(t, \sigma, \xi)=|\xi|^2t^m
\hat{T}(t,\sigma,\xi)$ in [8]. The main reason is: In the integrals (2.18)-(2.22), the variables $t$ and $x$ are the parameter variable
and  the integration variable respectively. This brings more troubles in treating the integrals (2.18)-(2.22)
of $\hat{K}(t, \sigma, \xi)$ than in treating the corresponding integrals of $\hat{K_0}(t, \sigma, \xi)$
due to the appearance of the integral variable factor $\si^m$ in $\hat{K}(t, \sigma, \xi)$.

\vskip .3cm
{\bf Lemma 2.5.} {\it Let $\hat{K}(t, \sigma, \xi)$ be defined in (2.16), then

(a) (2.18) and (2.19) hold.

(b) (2.20) holds for $\al=0$, namely, for any $h(\si) \in C^\infty_0(\R^+)$, one has
$$\sup_{\f r 2 \le |\xi|\le 2 r} \int_{0}^{+\infty} \Big|
\int_{0}^{+\infty} \hat{K}(t, \sigma, \xi) h(\sigma) d\sigma \Big| dt \le C \int_0^{+\infty}
|h(\sigma)| d\sigma.\eqno{(2.24)}
$$
}

\vskip .1cm

{\bf Proof.} (a)  Noting that
\begin{align*} \int_0^{+\infty} |\hat{K}(t,\sigma,\xi)| d\sigma =&
\int_0^{+\infty} |\xi|^2 \sigma^m \Big( \int_0^{\min(t,\sigma)}
\f {
\la(ts) \la(\sigma s)} { \lambda^2(ys)} dy\Big) d\sigma \\
=& \la(ts) \int_0^{ts} \f {dy} {\la^2(y)} \int_y^{+\infty} \si^m
\la(\si) d\si \qquad (\text{due to $|\xi|^2=s^{m+2}$}) \\
=& \la(ts) \int_0^{ts} \f {dy} {\la^2(y)} \int_y^{+\infty}
\la''(\si) d\si  \qquad \quad ( \text{due to} \ \la''(\si)=\si^m \la(\si)) \\
=& 1-\la(ts)\le C,\qquad\qquad \qquad\qquad  \qquad \qquad \qquad \qquad  ( \text{by (2.2)})
\end{align*}
then (2.18) holds.

In addition, one has
\begin{align*}\int_0^{+\infty} |\hat{K}(t,\si,\xi)| dt
=& \int_0^{+\infty} |\xi|^2 \si^m \Big(\int_0^{\min(t,\sigma)} \f {
\la(ts) \la(\sigma
s)} { \lambda^2(ys)} dy \Big) dt \notag\\
=& (\si s)^m  \int_0^{\si s} dy \int_y^{+\infty} \f {\la(\si s)
\la(t)} {\la^2(y)} dt.
\end{align*}
Without loss of generality, we may assume $\si s \ge 2$ (otherwise, the uniform bound of
$\int_0^{+\infty} |\hat{K}(t,\si,\xi)| dt$ can be easily obtained since $
\int_0^{+\infty} |\hat{K}(t,\si,\xi)| dt \le 2^m \la(\si s) \int_0^2 \f{dy}{\la^2(y)} \int_0^{+\infty} \la(t)
dt \le C$ holds). In this case, we have
\begin{align*} \int_0^{+\infty} |\hat{K}(t,\si,\xi)| dt=& (\si
s)^m\Big(\int_0^1 + \int_{1}^{\f {\si s} 2} + \int_{\f {\si s}
2}^{\si s} \Big) dy \int_y^{+\infty} \f { \la(t) \la(\sigma s)} {
\lambda^2(y)} dt\\ \equiv & L_1+L_2+L_3.\tag{2.25}
\end{align*}
Next we treat each $L_i$ $(i=1,2,3)$ in (2.25). From (2.1) and (2.2), one has
$$L_1\le C\la(\si s) (\si s)^m\la^{-2}(1) \int_0^{+\infty} (1+t)^{-2} dt\le
C.  \eqno{(2.26)}$$
The boundedness of $L_2$ and $L_3$ can  be obtained by applying
(2.3). Indeed,
\begin{align*}   L_2 \le & (\si s)^m\int_{1}^{\f {\si s} 2} dy
\int_y^{+\infty} \big( \f t y \big)^{\f 1 2}  \big( \f {\si s} y
\big)^{\f 1 2} exp\Big( \f 4 {m+2} y^{\f {m+2} 2 } -\f 2 {m+2} t^{\f
{m+2} 2 }-\f 2 {m+2} (\si s)^{\f
{m+2} 2 }  \Big) dt \\
\le & (\si s)^{m+\f 1 2} \ exp\Big( -\f 2 {m+2} (\si s)^{\f {m+2} 2
}\Big) \int_{1}^{\f {\si s} 2} y^{-1}exp\Big( \f 4 {m+2} y^{\f {m+2}
2 }\Big)  dy \\
& \quad \times\int_y^{+\infty} t^{\f {1-m} 2} exp\Big(
-\f 2 {m+2} t^{\f {m+2} 2 }\Big) d(t^{\f {m+2} 2}) \notag\\
\le & (\si s)^{m+\f 1 2} \ exp\Big( -\f 2 {m+2} (\si s)^{\f {m+2} 2
}\Big) \int_{1}^{\f {\si s} 2} y^{-\f {m+1}{2}}exp\Big( \f 2 {m+2}
y^{\f {m+2} 2 }\Big) dy \\
\le & (\si s)^{m+\f 1 2} \ exp\Big( -\f 2 {m+2} (\si s)^{\f {m+2} 2
}\Big) \ exp\Big( \f 2 {m+2} (\f {\si s} {2})^{\f {m+2} 2 }\Big) \notag\\
\le & exp\Big( -C (\si s)^{\f {m+2} 2 }\Big) \\
\le & C \tag{2.27}
\end{align*}
and
\begin{align*}
L_3 &
\le  (\si s)^{m+\f 1 2} \ exp\Big( -\f 2 {m+2} (\si s)^{\f {m+2} 2
}\Big) \int_{\f {\si s} 2}^{\si s} y^{-1}exp\Big( \f 4 {m+2} y^{\f
{m+2} 2 }\Big) dy \\
& \quad \times \int_y^{+\infty} t^{\f {1-m} 2} exp\Big(
-\f 2 {m+2} t^{\f {m+2} 2 }\Big) d(t^{\f {m+2} 2}) \\
& \le  (\si s)^{m+\f 1 2} \ exp\Big( -\f 2 {m+2} (\si s)^{\f {m+2} 2
}\Big) \int_{\f {\si s} 2}^{ \si s} y^{-m- \frac{ 1}{ 2}}exp\Big( \f
2 {m+2}
y^{\f {m+2} 2 }\Big) d(y^{ \f {m+2} 2}) \\
&\le (\si s)^{m+\f 1 2} (\f{\si s} 2)^{-m-\f 1 2}\le C.  \tag {2.28}
\end{align*}
Substituting (2.26)-(2.28) into (2.25) yields (2.19).
(b) It follows from (a) that
$$
\int_0^{+\infty} \Big| \int_0^{+\infty} \hat{K}(t,\si,\xi)
h(\si) d\si \Big| dt \le  \int_0^{+\infty} |h(\si)| \Big(
\int_0^{+\infty}   \hat{K}(t,\si,\xi)  dt \Big) d\si
\le C \int_0^{+\infty} |h(\si)| d\si.
$$
Thus (2.20) holds for $\alpha=0.$ \qquad  \qquad \hfill
 $\square$
\vskip .3cm

{\bf Lemma 2.6.} {\it  Let $\hat{K}(t, \sigma, \xi)$ be defined in (2.16), then
(2.21) holds for $\al=0$, i.e.,
$$\sup_{\f r 2 \le |\xi|\le 2 r} \int_{\mathcal C \Delta(a, 2^{q+1}b)} \Big|
\int_0^{+\infty} \hat{K}(t, \sigma, \xi) h(\sigma) d\sigma \Big| dt \le
C_q P_1(ar^{\f{2}{m+2}}, br^{\f{2}{m+2}})\int_0^{+\infty} |h(\sigma)| d\sigma,\eqno{(2.29)}
$$
where $h(\si)\in C_0^{\infty}(\Bbb R^+)$ is supported in $\Delta=(a-b, a+b)$, and $P_1(a, b)=
exp\Big(-C a^{m/2}b\Big)$.}

\vskip .1cm

{\bf Proof.} Obviously, it suffices to establish (2.29) for $q=0$ if the constant $C_q$
can be shown to be independent of $q$. Notice that \begin{equation*}
\mathcal C \Delta(a, 2b)=\left\{
\begin{aligned}
& (a+2b, +\infty),\qquad \qquad \qquad \text{if $a\le 2b$}\\
&(a+2b, +\infty)\cup (0, a-2b),\quad\text{if $a>2b$},\\
\end{aligned}
\right.
\end{equation*} then in order to estimate the integral $\int_{\mathcal C \Delta(a, 2^{q+1}b)}$
in (2.29), we require to deal with the integrals $\int_{a+2b}^{+\infty}$
and  $\int_{0}^{a-2b}$ separately.

It follows from a direct computation that
\begin{align*} & \int_{a+2b}^{+\infty} \Big| \int_0^{+\infty}
\hat{K}(t,\si,\xi) h(\si) d\si \Big|  dt\\
\le& \int_{a+2b}^{+\infty} dt \Big( \int_{a-b}^{a+b} |h(\si)| \
|\xi|^2 \sigma^m d\si\big( \int_0^{\min(t,\sigma)} \f { \la(ts)
\la(\sigma s)} { \lambda^2(ys)} dy \big) \Big) \\
=& \int_{a-b}^{a+b} |h(\si)| \  (\sigma s)^m \la(\si s) d\si
\Big(\int_0^{\si s} \f {dy}{\la^2(y)} \int_{(a+2b)s}^{+\infty}
\la(t) dt\ \Big) \\
 \le & \int_{a-b}^{a+b} |h(\si)| \  (\sigma s)^{m+1}d\si \Big(
\int_{(a+2b)s}^{+\infty} \f {\la(t)}{\la(\si s)} dt \Big) \qquad
\qquad  \text{(by (2.1))}\\
\le & \int_{a-b}^{a+b} |h(\si)| \  (\sigma s)^{m+1} d \si\Big(
\int_{(a+2b)s}^{+\infty} \big(\f t {\si s} \big)^{\f 1 2} exp\big(
\f 2 {m+2} (\si s)^{\f {m+2} 2 } -\f 2 {m+2} t^{\f {m+2} 2 }\big) dt
\Big) \quad (\text{ by  (2.3)}) \\
\le & \int_{a-b}^{a+b} |h(\si)| \ (\sigma s)^{m+\f 1 2} \ exp\big(
\f 2 {m+2} (\si s)^{\f {m+2} 2 }\big) d\si \int_{(a+2b)s}^{+\infty}
t^{\f {1-m} 2}  \ exp\big( -\f 2 {m+2} t^{\f {m+2} 2 }\big) d(t^{\f
{m+2} 2}) \\
\le & \int_{a-b}^{a+b} |h(\si)| \ (\sigma s)^{m+\f 1 2} \
\big((a+2b)s \big)^{\f {1-m} 2} \ exp\big(\f 2 {m+2} (\si s)^{\f
{m+2} 2 }-\f 2 {m+2} ((a+2b)s)^{\f {m+2} 2 }\big) d\si \\
\le & C \big((a+b)s \big)^{\f {m} 2 +1} \ exp\big(-C (as)^{\f m 2}
(bs) \big) \int_0^{+\infty} |h(\si)|d\si \\
\le & C exp\Big(-C (as)^{\f m 2} (bs) \Big) \int_0^{+\infty} |h(\si)|d\si, \tag{2.30}
\end{align*}
where the positive constant $C$ is independent of $a$ and $b$.

If $a-2b \ge 0,$ then
\begin{align*} &\int_0^{a-2b} dt \Big| \int_0^{+\infty} \hat{K}(t, \sigma,
\xi) h(\si) d\si  \Big| \\ \le & \int_{a-b}^{a+b} |h(\si)| \ (\si
s)^m \la(\si s) d\si \Big( \int_0^{(a-2b)s} \la(t) dt\int_0^t \f {dy} {\la^2(y)} \Big) \\
\le & \int_{a-b}^{a+b} |h(\si)|  (\si s)^m \la(\si s) d\si
\int_0^{(a-2b)s} \f {t} {\la(t)} dt \qquad  \text{(by (2.1))}\\
=& \int_0^{(a-2b)s} t dt \Big(\int_{a-b}^{a+b} |h(\si)|  (\si s)^m
\f{\la(\si s)} {\la(t)} d\si \Big)   \\
\le & \int_0^{(a-2b)s} t dt \Big(\int_{a-b}^{a+b} |h(\si)|  (\si
s)^m \big(\f{\si s}{t}\big)^{1/2} exp\big( \f 2 {m+2}t^{\f {m+2} 2 }
-\f 2 {m+2} (\si s)^{\f {m+2} 2 } \big) d\si \Big) \quad ( \text{ by (2.3)} )\\
\le & C \big( (a-2b)s\big)^{\f 3 2}\big( (a+b)s\big)^{m+\f 1 2}  \ exp\big( \f 2
{m+2}((a-2b)s)^{\f {m+2} 2 } \big) \
exp\big(- \f 2 {m+2}((a-b)s)^{\f {m+2} 2 } \big)\\
& \quad \times\int_0^{+\infty}
|h(\si)| d\si\\
\le & C exp\big(-C (as)^{\f m 2}(bs) \big)  \int_0^{+\infty} |h(\si)| d\si. \tag{2.31}
\end{align*}
From (2.30) and (2.31), we see that
$$ \int_{C \Delta(a, 2^{q+1}b)} \Big|
\int_0^{+\infty} \hat{K}(t, \sigma, \xi) h(\sigma) d\sigma \Big| dt
\le C P_1(as,bs)\int_0^{+\infty} |h(\sigma)| d\sigma,\eqno{(2.32)}$$ where $P_1(a, b)=
exp\Big(-C a^{m/2}b)\Big)$. Thus, (2.29) obviously holds due to (2.32).
\hfill $\square$

\vskip .3cm

{\bf Lemma 2.7.} {\it  Let $\hat{K}(t, \sigma, \xi)$ be defined in (2.16), then
(2.22) holds for $\al=0$, i.e.,
$$\sup_{\f r 2 \le |\xi|\le 2 r} \int_0^{+\infty} \Big|
\int_0^{+\infty} \hat{K}(t, \sigma, \xi) h(\sigma) d\sigma \Big| dt \le
CP_2(ar^{\f{2}{m+2}}, br^{\f{2}{m+2}})\int_0^{+\infty} |h(\sigma)| d\sigma,\eqno{(2.33)}
$$
where $h(\si)\in C_0^{\infty}(\Bbb R^+)$, $supp h(\si)\subset\Delta=(a-b, a+b)$, $\int_0^{+\infty} h(\si)d\si=0$, and
$P_2(a,b)=b a^{\f m 2} + b^2 a^m +b^3 a^{\f
{3m} 2}.$}

\vskip .1cm
{\bf Proof.} We will divide the proof of (2.33) into the following three parts:

$$\biggl(\int_{a+b}^{+\infty}+\int_0^{a-b}+\int_{a-b}^{a+b}\biggr)\Big| \int_0^{+\infty} \hat{K}(t,
\sigma, \xi) h( \si) d\si \Big| dt.$$

\vskip 0.2 true cm

{\bf Part 1. $ \int_{a+b}^{+\infty} \Big| \int_0^{+\infty} \hat{K}(t,
\sigma, \xi) h( \si) d\si \Big| dt\le C bs (as)^{\f m 2} \int_0^{+\infty}
|h(\si)| d\si$ holds}

\vskip 0.2 true cm

By a direct computation, we have
\begin{align*} & \int_{a+b}^{+\infty} \Big| \int_0^{+\infty} \hat{K}(t,
\sigma, \xi) h( \si) d\si \Big| dt \\
=& \int_{a+b}^{+\infty} \Big| \int_{a-b}^{a+b} \big( \hat{K}(t,
\sigma, \xi)- \hat{K}(t, a+b, \xi)  \big) h(\si) d\si   \Big| dt
\qquad (\text{by }\int h(\si)d\si =0)
\\
=& \int_{a+b}^{+\infty} \la(ts) dt \Big| \int_{a-b}^{a+b} h(\si)
|\xi|^2 \Big( \si^m \la(\sigma s)\int_0^{\sigma} \f { dy } {
\lambda^2(ys)} - (a+b)^m \la((a+b) s)\int_0^{a+b} \f { dy } {
\lambda^2(ys)} \Big)d\si  \Big|
\\
\le & \int_{(a+b)s}^{+\infty} \la(t) dt\int_{a-b}^{a+b} |h(\si)|
(\si s)^m \Big( \la(\si s) -\la((a+b)s) \Big) d\si \int_0^{\si s} \f
{dy}{\la^2(y)} \\
& \   + \int_{(a+b)s}^{+\infty} \la(t) dt \int_{a-b}^{a+b} |h(\si)|
\Big(\big((a+b)s\big)^m - \big(\si s\big)^m \Big) \la((a+b)s)
 d\si \int_0^{\si s} \f {dy}{\la^2(y)} \\
 & \ +\int_{(a+b)s}^{+\infty} \la(t) dt\int_{a-b}^{a+b} |h(\si)|
 \big((a+b)s \big)^m \la((a+b)s) d\si \int_{\si s}^{(a+b)s} \f {dy}
 {\la^2(y)} \\
 \equiv& I_1 +I_2 +I_3.\tag{2.34}
\end{align*}
Next we treat each $I_i$ ($i=1, 2, 3$) in (2.34). For $I_1$, one has
\begin{align*}  I_1=& \int_{a-b}^{a+b} |h(\si)| (\si s)^m \la(\si s) d\si
\Big(\int_0^{\si s} \f {\la(\si s) -\la((a+b)s)}{ \la^2(y)}dy \Big)
\int_{(a+b)s}^{+\infty} \f {\la(t)} {\la(\si s)} dt\\
\le & C bs\int_{a-b}^{a+b} |h(\si)| (\si s)^m \la(\si s) d\si
\Big(\int_0^{\si s}\f{-\la'(y)}{\la^2(y)} dy \Big)
\int_{(a+b)s}^{+\infty} \big( \f t {\si s} \big)^{\f 1 2}  \\
& \quad \times exp\big(
\f 2 {m+2} (\si s)^{\f {m+2} 2 } -\f 2 {m+2} t^{\f {m+2} 2 }  \big)
dt \qquad  ( \text{by (2.1)  and  (2.3)}) \\
=& C bs\int_{a-b}^{a+b} |h(\si)| (\si s)^{m-\f 1 2}  \
exp\big( \f 2 {m+2} (\si s)^{\f {m+2} 2 } \big) \big(1- \la(\si s)
\big) d\si \\
& \quad \times \int_{(a+b)s}^{+\infty} t^{\f {1-m}{2}}exp\big( -\f 2
{m+2} t^{\f
{m+2} 2 } \big) d(t^{\f {m+2} 2}) \\
\le & C bs \big( (a+b)s \big)^{\f {1-m} 2} exp\big( -\f 2 {m+2}
((a+b)s)^{\f {m+2} 2 } \big) \int_{a-b}^{a+b} |h(\si)| (\si s)^{m-\f 1 2}  \\
& \quad \times exp\big(
\f 2 {m+2} (\si s)^{\f {m+2} 2 } \big) \big(1- \la(\si s) \big) d\si
\\
\le & C bs((a+b)s)^{\f m 2} \int_0^{+\infty} |h(\si)| \big(1- \la(\si s) \big)d\si \\
\le & C bs(as)^{\f m 2} \int_0^{+\infty} |h(\si)|d\si.  \qquad ( \text{by
(2.2)}) \tag {2.35}
\end{align*}
Similarly, applying Lemma 2.2 yields
\begin{align*}  I_2 \le & C bs\big((a+b)s\big)^{m-1} \la((a+b)s)\int_{a-b}^{a+b}
|h(\si)| d\si \Big(\int_0^{\si s} \f {\la(\si s)}{\la^2(y)} dy \Big)
\Big( \int_{(a+b)s}^{+\infty}  \f {\la (t)}{\la(\si s)}dt\Big) \\
\le & C bs\big((a+b)s\big)^{m-1} \la((a+b)s)\int_{a-b}^{a+b}
|h(\si)| \f 1{\la(\si s)}\Big(\int_0^{\si s} \f {\la(\si s)}{\la(y)}
dy \Big) \Big( \int_{(a+b)s}^{+\infty}  \f {\la (t)}{\la(\si
s)}dt\Big) \\
\le & C bs\big((a+b)s\big)^{m-1}\int_{a-b}^{a+b} |h(\si)| d\si
\Big(\int_0^{\si s}  \big( \f {\si s} {y} \big)^{\f 1 2} exp\big( \f
2 {m+2} y^{\f {m+2} 2 } -\f 2 {m+2} (\si s)^{\f {m+2} 2 }  \big)dy
\Big) \\
& \ \times \Big( \int_{(a+b)s}^{+\infty}\big( \f t {\si s} \big)^{\f
1 2}  exp\big( \f 2 {m+2} (\si s)^{\f {m+2} 2 } -\f 2 {m+2} t^{\f
{m+2} 2 }  \big) dt \Big) \qquad (\text{by (2.1) and (2.3) }) \\
\le & C bs\big((a+b)s\big)^{\f{m-1}{2}} \int_{a-b}^{a+b} |h(\si)|
exp\big( \f 2 {m+2} (\si s)^{\f {m+2} 2 } \big)  d\si \Big(
\int_0^{\si s} y^{-1/2} dy\Big)\\
& \quad \times \Big(
\int_{(a+b)s}^{+\infty} exp\big(
-\f 2 {m+2} t^{\f {m+2} 2 } \big) dt \Big) \\
\le & C bs (a s)^{\f m 2 } \int_0^{+\infty} |h(\si)| d\si \tag {2.36}
\end{align*}
and
\begin{align*}
 I_3 =& \big((a+b)s \big)^m  \int_{a-b}^{a+b}
|h(\si)| d\si \int_{\si s}^{(a+b)s} \f {\la^2((a+b)s)}
 {\la^2(y)} dy \int_{(a+b)s}^{+\infty}  \f {\la(t)}{\la((a+b)s)} dt\\
 \le & C bs\big((a+b)s \big)^m \int_{a-b}^{a+b}
|h(\si)| d\si \int_{(a+b)s}^{+\infty} \big( \f t {(a+b) s} \big)^{\f
1 2} \\
& \quad \times exp\big( \f 2 {m+2} ((a+b) s)^{\f {m+2} 2 } -\f 2 {m+2} t^{\f
{m+2} 2 }  \big)dt \quad (\text{by (2.3)})\\
 =& C bs
\big((a+b)s \big)^{m- \f 1 2} exp\big( \f 2 {m+2} ((a+b) s)^{\f
{m+2} 2 } \big)\int_{a-b}^{a+b} |h(\si)| d\si\int_{(a+b)s}^{+\infty}
t^{\f {1-m} 2} \\
& \quad \times exp\big( -\f 2 {m+2} t^{\f {m+2} 2 }  \big) d(t^{\f
{m+1} 2}) \\
\le & C bs (as)^{\f m 2} \int_0^{+\infty} |h(\si)| d\si.  \tag {2.37}
\end{align*}

Substituting (2.35)-(2.37) into (2.34) yields the conclusion in Part 1.

\vskip 0.2 true cm

{\bf Part 2. $\int_0^{a-b} \Big| \int_0^{+\infty} \hat{K}(t, \sigma,
\xi) h(\si) d\si \Big| dt \le C bs (as)^{ \f m 2}\int_0^{+\infty} |h(\si)|
d\si$ holds}

\vskip 0.2 true cm
Note that
\begin{align*} & \int_0^{a-b} \Big| \int_0^{+\infty} \hat{K}(t, \sigma,
\xi) h(\si) d\si \Big| dt \\
=& \int_0^{a-b} \Big| \int_0^{+\infty} \Big( \hat{K}(t, \sigma,
\xi)- \hat{K}(t, a+b, \xi)\Big) h(\si) d\si \Big| dt \qquad  (\text{by} \int_0^{+\infty} h(\si)d\si =0 )\\
=& \int_0^{(a-b)s} \la(t) dt\Big| \int_{a-b}^{a+b} h(\si) \Big( (\si
s)^m  \la(\sigma s) \int_0^{t} \f {dy} { \lambda^2(y)}  -
((a+b)s)^m \la((a+b) s)\int_0^{t} \f {dy} { \lambda^2(y)}  \Big)
d\si \Big|
\\
\le &   \int_{a-b}^{a+b} |h(\si)| (\si s)^m \Big( \la(\si s)
-\la((a+b)s) \Big)d\si \int_0^{(a-b)s} \la(t) dt \int_0^t \f {dy } {
\lambda^2(y)} \\
& \ +  \la((a+b)s) \int_{a-b}^{a+b} |h(\si)| \Big((\si s)^m
-((a+b)s)^m \Big)d\si  \int_0^{(a-b)s} \la(t) dt \int_0^t \f {dy } {
\lambda^2(y)} \\
\le & C bs\int_{a-b}^{a+b} |h(\si)| (\si s)^m |\la'(\si s)| d\si
\int_0^{(a-b)s} \la(t) dt\int_0^t\f {dy} { \lambda^2(y)}\\
& \ +C bs ((a+b)s)^{m-1} \la((a+b)s) \int_{a-b}^{a+b} |h(\si)|
d\si \int_0^{(a-b)s} \la(t) dt \int_0^t \f {dy} { \lambda^2(y)}
\quad  (\text{by
(2.1)})\\
\equiv& II_1 +II_2.\tag{2.38}
\end{align*}
To estimate the term $II_1$ in (2.38), we will consider the following three cases

\vskip 0.2 true cm

{\bf Case 1. $(a-b)s \ge 2$ }

\vskip 0.2 true cm

One now has
\begin{align*}
II_1&\le C bs\int_{a-b}^{a+b} |h(\si)| (\si s)^m |\la'(\si
s)| d\si \biggl(\int_0^1 \la(t) dt \int_0^t \f {dy} {\la^2(y)}+\int_1^{\f{(a-b)s} 2} \la(t) dt \int_0^t \f {dy}
{\la^2(y)}\\
&+ \int_{\f{(a-b)s} 2}^{(a-b)s} \la(t) dt \int_0^1
\f {dy} {\la^2(y)}+ \int_{\f{(a-b)s} 2}^{(a-b)s} \la(t) dt\int_1^{\f
{(a-b)s} 2} \f {dy} {\la^2(y)}+\int_{\f{(a-b)s} 2}^{(a-b)s} \la(t) dt \int_{\f
{(a-b)s} 2}^{t} \f {dy} {\la^2(y)}\biggr)\\
&\equiv II_1^{(1)}+II_1^{(2)}+II_1^{(3)}+II_1^{(4)}+II_1^{(5)}.\tag{2.39}
\end{align*}
We only estimate $II_1^{(1)}$, $II_1^{(4)}$ and $II_1^{(5)}$ since
the remaining terms can be  treated similarly.

It follows from Lemma 2.2 that
\begin{align*}  II_1^{(1)} \le & C (bs)\int_{a-b}^{a+b} |h(\si)| (\si s)^m
|\la'(\si s)| d\si \int_0^1  \la(t) \f t {\la^2(t)} dt \quad (  \text{by (2.1)})\\
 \le &  C \la(1)^{-1} \ bs \int_{a-b}^{a+b} |h(\si)| (\si s)^m
|\la'(\si s)| d\si \qquad (  \text{by (2.1)})\\
\le & C bs (as)^{\f m 2} \int_0^{+\infty} |h(\si)|d\si\quad ( \text{by
(2.2)}) \tag {2.40}
\end{align*}
and
\begin{align*}  II_1^{(4)}\le & C bs\int_{a-b}^{a+b} |h(\si)| (\si s)^m \f
{|\la'(\si s)|}{\la(\si s)} d\si \int_{\f{(a-b)s} 2}^{(a-b)s}
 dt \int_1^{\f {(a-b)s} 2} \f { \la(\si s)  \la(t) } {\la^2(y)} dy \\
 \le & C bs\int_{a-b}^{a+b} |h(\si)|  (\si s)^{ \f {3m} 2}  d\si \int_{\f{(a-b)s} 2}^{(a-b)s}
 dt \\
 & \ \times \int_1^{\f {(a-b)s} 2} \big( \f t y \big)^{\f 1 2}  \big( \f {\si s} y \big)^{\f 1 2}
exp\Big( \f 4 {m+2} y^{\f {m+2} 2 } -\f 2 {m+2} t^{\f {m+2} 2 }-\f 2
{m+2} (\si s)^{\f {m+2} 2 }  \Big) dy \\
& \hskip 6cm (\text{  by (2.4) since $(a-b)s \ge 2$ and by (2.3)})\\
\le & C bs\int_{a-b}^{a+b} |h(\si)|(\si s)^{ \f {3m+1}
2}exp\Big(-\f 2 {m+2} (\si s)^{\f {m+2} 2 }  \Big)
\big(\f{(a-b)s} 2\big)^{\f {3} 2} \\
& \quad \times exp\Big(\f 2 {m+2}
\big(\f {(a-b) s)}{2 }\big)^{\f {m+2} 2 }  \Big)d\si \\
\le & C bs (as)^{\f m 2} \int_0^{+\infty} |h(\si)| d\si. \tag {2.41}
\end{align*}
In addition, by (2.3) and (2.4), we arrive at
\begin{align*}
II_1^{(5)}= & C bs\int_{a-b}^{a+b} |h(\si)| (\si s)^m \f
{|\la'(\si s)|}{\la(\si s)} d\si \int_{\f{(a-b)s} 2}^{(a-b)s} dy
\int_{\f {(a-b)s} 2}^{y} \f {\la(\si s) \la(y)} {\la^2(t)} dt \\
\le & C bs\int_{a-b}^{a+b} |h(\si)|(\si s)^{ \f {3m+1}
2}exp\Big(-\f 2 {m+2} (\si s)^{\f {m+2} 2 }  \Big)
d\si\\
& \ \times \int_{\f{(a-b)s} 2}^{(a-b)s}t^{\f {1-m} 2}exp\Big(-\f 2
{m+2} t^{\f {m+2} 2 }  \Big) d(t^{\f {m+1} 2}) \int_{\f {(a-b)s}
2}^{t} y^{-1-\f {m} 2}exp\Big(\f 4 {m+2} y^{\f {m+2} 2 }  \Big)
d(y^{\f {m+2} 2}) \\
\le & C bs\big((a+b)s\big)^{\f m 2} \int_{a-b}^{a+b} |h(\si)|(\si
s)^{ m + \f 1 2}exp\Big(-\f 2 {m+2} (\si s)^{\f {m+2} 2 } \Big) d\si
\\& \ \times \big((a-b) s \big)^{-( m + \f 1 2 )}exp\Big(\f 2 {m+2} \big((a-b) s
\big)^{\f {m+2} 2 } \Big). \tag {2.42}
\end{align*}
Notice that the function $\eta(z)=z^{m+ \f 1 2} exp\Big(-\f 2
{m+2} z^{\f {m+2} 2 } \Big)$ is strictly decreasing for $z\ge  2,$
then (2.42) can be dominated by
$$
C bs \big((a+b)s\big)^{\f m 2} \int_0^{+\infty} |h(\si)| d\si \le C bs (as)^{\f m 2} \int_0^{+\infty} |h(\si)|d\si.
$$
Collecting all the analysis above yields
$$II_1 \le C bs (as)^{\f m 2} \int_0^{+\infty} |h(\si)|
d\si\qquad \text{for $(a-b)s \ge 2$}.\eqno{(2.43)}$$

{\bf Case 2. $(a-b)s \le 1$}

\vskip 0.2 true cm

In this case, we have $$II_1\le II_1^{(1)}\le C bs (as)^{\f m 2} \int_0^{+\infty}|h(\si)|
d\si.\eqno{(2.44)}$$

{\bf Case 3. $1 \le (a-b)s \le 2$}

\vskip 0.2 true cm

In this case, we have $$II_1 \le II_1^{(1)} +
II_1^{(6)},$$
where $II_1^{(6)}=C bs \int_{a-b}^{a+b} |h(\si)| (\si s)^m
|\la'(\si s)| d\si \int_1^2 \la(t) dt \int_0^t \f {dy} {\la^2(y)}$.

As in Case 2, one has $II_1^{(1)}\le C bs (as)^{\f m 2} \int_{\Bbb R} |h(\si)|
d\si$. In addition,
\begin{align*} II_1^{(6)}\le & C bs\int_{a-b}^{a+b} |h(\si)| (\si s)^m
|\la'(\si s)| d\si \int_1^2 \f t {\la(t)} dt  \quad (\text{by
(2.1)}) \\
\le & C bs (as)^{\f m 2} \int_{\Bbb R} |h(\si)|
  d\si.
\end{align*}
Therefore,
$$ II_1\le C bs (as)^{\f m 2} \int_0^{+\infty} |h(\si)|
d\si\quad\text{for $1 \le (a-b)s \le 2$}.\eqno{(2.45)}$$
Combining (2.43)-(2.45) yields
$$II_1\le C bs (as)^{\f m 2} \int_0^{+\infty} |h(\si)|
 d\si.\eqno{(2.46)}$$
Next we estimate $II_2$ in (2.38), which will be divided into the following two cases.

\vskip 0.2 true cm

{\bf Case A. $(a-b)s \ge 1 $}

\vskip 0.2 true cm

We have
\begin{align*}
II_2 \le &  C bs \big((a+b)s \big)^{m-1} \la((a+b)s)
\int_{a-b}^{a+b} |h(\si)| d\si\\
& \ \ \times\biggl(\int_0^1  \la(t) dt \int_0^t \f {dy}
{\la^2(y)} +\int_1^{(a-b)s}  \la(t) dt \int_0^1
\f {dy} {\la^2(y)}+\int_1^{(a-b)s}  \la(t) dt \int_1^t
\f {dy} {\la^2(y)}\biggr)\\
\equiv & II_2^{(1)} + II_2^{(2)} + II_2^{(3)}.\tag{2.47}
\end{align*}
Noting the function $\la(t)$ is decreasing, then by (2.1)-(2.3)
\begin{align*}  II_2^{(1)}\le &C bs \big((a+b)s \big)^{m-1}  \int_{a-b}^{a+b} |h(\si)| d\si
\int_0^1 t \big( \f  {(a+b) s} t \big)^{1/2} \\
& \quad \times exp\big( \f 2 {m+2}
t^{\f {m+2} 2 }- \f 2 {m+2} \big((a+b) s \big)^{\f {m+2} 2 } \big)
dt \\
\le & C bs \big((a+b)s \big)^{m-\f 1 2} exp\big( -\f 2 {m+2}
\big((a+b)s\big)^{\f {m+2} 2 } \big) \int_{a-b}^{a+b} |h(\si)| d\si
\\
& \quad \times  \int_0^1 t^{1/2} exp\big( \f 2 {m+2} t^{\f {m+2} 2 } \big) dt\\
\le & C bs (as)^{\f m 2} \int_0^{+\infty} |h(\si)| d\si \tag {2.48}
\end{align*}
and
\begin{align*}   II_2^{(2)}\le &  C \f {bs}{\la^2(1)} \big((a+b)s
\big)^{m-1} \la((a+b)s) \int_{a-b}^{a+b} |h(\si)| d\si
\int_1^{+\infty} \la(t) dt \\
 \le & C bs\big((a+b)s \big)^{\f m 2}
 \Big((a+b)s \Big)^{\f m 2 -1}\la((a+b)s)\int_0^{+\infty} |h(\si)| d\si\\
\le & C bs (as)^{\f m 2 } \int_0^{+\infty} |h(\si)|
d\si. \tag {2.49}
\end{align*}
For the term $ II_2^{(3)}$, one has from (2.3)
\begin{align*} II_2^{(3)}= &  C bs\big((a+b)s \big)^{m-1}
\int_{a-b}^{a+b} |h(\si)| d\si  \int_1^{(a-b)s}   dt \int_1^t \f
{  \la(t) \la((a+b)s)   } {\la^2(y)} dy \\
\le & C bs\big((a+b)s \big)^{m-\f 1 2} exp\Big(-\f 2 {m+2}
\big((a+b) s\big)^{\f {m+2} 2 }  \Big)\int_{a-b}^{a+b} |h(\si)|
d\si \\
& \quad \times \int_1^{(a-b)s} t^{1/2} exp\Big(\f 2 {m+2} t^{\f {m+2} 2
}  \Big)  dt \\
\le & C bs \big((a+b)s \big)^{m-\f 1 2}\big((a-b)s \big)^{\f 3 2}\int_{a-b}^{a+b} |h(\si)|
d\si
\\
& \quad \times exp\Big(\f 2 {m+2} \big((a-b) s\big)^{\f {m+2} 2 }-\f 2 {m+2}
\big((a+b) s\big)^{\f {m+2} 2 } \Big) \\
\le & C bs (as)^{ \f m 2}\int_0^{+\infty} |h(\si)| d\si. \tag {2.50}
\end{align*}
From (2.48)-(2.50), we see that
$$II_2 \le C bs (as)^{ \f m 2}\int_0^{+\infty} |h(\si)| d\si \qquad \text{ for } (a-b)s \ge 1.
\eqno{(2.51)}$$

{\bf Case B. $(a-b)s \le 1 $}
\vskip 0.2 true cm

One easily obtains
$$II_2 \le II_2^{(1)} \le C bs (as)^{ \f m 2}\int_0^{+\infty} |h(\si)| d\si \qquad \text{ for } (a-b)s \le 1.
\eqno{(2.52)}$$

Therefore, it follows from (2.38), (2.46) and (2.51)-(2.52) that the conclusion in Part 2 holds.

\vskip 0.3 true cm

{\bf Part 3. $\int_{a-b}^{a+b} \Big| \int_0^{+\infty} \hat{K}(t, \sigma,
\xi) h(\si) d\si \Big| dt \le C\Big( bs(as)^{\f m 2}
+(bs)^2 (as)^m +(bs)^3 (as)^{\f {3m} 2} \Big)\int_0^{+\infty} |h(\si)|
d\si$ holds}

We have
\begin{align*} & \int_{a-b}^{a+b} \Big| \int_0^{+\infty} \hat{K}(t,
\sigma, \xi) h( \si) d\si \Big| dt \\
=& \int_{a-b}^{a+b} \Big| \int_{a-b}^{a+b} h(\si) |\xi|^2 \Big(
\si^m \int_0^{\min(t,\sigma)} \f { \la(ts) \la(\sigma s)} {
\lambda^2(ys)} dy -  (a+b)^m \int_0^{t} \f { \la(ts) \la((a+b) s)} {
\lambda^2(ys)} dy
\Big)d\si \Big| dt \\
\le & \int_{(a-b)s}^{(a+b)s} \la(t) dt \int_{a-b}^{a+b} |h(\si)| (\si
s)^m \Big( \la(\si s) -\la((a+b)s) \Big) d\si \Big(\int_0^{(a-b)s}+
\int_{(a-b)s}^{\min(t,\si s)} \Big)\f
{dy}{\la^2(y)} \\
&  + \int_{(a-b)s}^{(a+b)s} \la(t) dt \int_{a-b}^{a+b} |h(\si)|
\Big(\big((a+b)s\big)^m - \big(\si s\big)^m \Big) \la((a+b)s)
 d\si  \Big(\int_0^{(a-b)s}+
\int_{(a-b)s}^{\min(t,\si s)} \Big) \f {dy}{\la^2(y)} \\
 & +\int_{(a-b)s}^{(a+b)s} \la(t) dt \int_{a-b}^{a+b} |h(\si)|
 \big((a+b)s \big)^m \la((a+b)s) d\si \int_{\min(t,\si s)}^{t} \f {dy}
 {\la^2(y)} \\
 \equiv & III_1 + III_2 +III_3 +III_4 +III_5.\tag{2.53}
\end{align*}

We now treat each $III_i$ in (2.53) as follows.

\vskip 0.3 true cm
{\bf Step 1. Estimate of $III_1$}

\vskip 0.2 true cm

By applying Lemma 2.2, one obtains
\begin{align*} III_1=& \int_{a-b}^{a+b} |h(\si)| (\si s)^m  d\si  \int_0^{(a-b)s} \f
{\la(\si s) -\la((a+b)s)}{\la^2(y)} dy \int_{(a-b)s}^{(a+b)s} \la(t)
dt\\
\le & C \Big((a+b)s \Big)^m \int_{a-b}^{a+b} |h(\si)| d\si \
\Big(bs \int_0^{(a-b)s} \f {-\la'(y)}{\la^2(y)} dy \Big) \ \Big(
\la((a-b)s)\int_{(a-b)s}^{(a+b)s} dt \Big)\\   
\le & C (bs)^2 (as)^m  \Big(1-\la((a-b)s) \Big) \int_0^{+\infty} |h(\si)| d\si
\\
\le & C  (bs)^2 (as)^m \int_0^{+\infty} |h(\si)| d\si. \tag {2.54}
\end{align*}

{\bf Step 2. Estimate of $III_2$}
\begin{align*} III_2
\le & C bs  \int_{a-b}^{a+b} |h(\si)| (\si s)^m |\la'(\si s)| d\si
\int_{(a-b)s}^{\si s} \la(t) dt
\int_{(a-b)s}^t \f {dy}{\la^2(y)} \\
& \ +  C bs  \int_{a-b}^{a+b} |h(\si)| (\si s)^m |\la'(\si s)|
d\si \int_{\si s}^{(a+b)s} \la(t) dt \int_{(a-b)s}^{\si s} \f
{dy}{\la^2(y)}\\
 \equiv & III_2^{(1)} + III_2^{(2)}.\tag{2.55}
\end{align*}
If $(a-b)s \ge 1,$ then
\begin{align*} III_2^{(1)} \le  & C  (bs)^2  \int_{a-b}^{a+b} |h(\si)| (\si
s)^m |\la'(\si s)| d\si \int_{(a-b)s}^{\si s} \f 1 {\la(y)} dy \qquad  (  \text{by (2.1)})\\
\le & C (bs)^3 \int_{a-b}^{a+b} |h(\si)| (\si s)^m \f { |\la'(\si
s)| }{\la(\si s)} d\si \quad  (  \text{by (2.1)})
\\
\le & C (bs)^3 \int_{a-b}^{a+b} |h(\si)| (\si s)^{ \f {3m} 2}  d\si
\qquad \qquad (  \text{by (2.4)})
\\ \le & C (bs)^3  (as)^{ \f {3 m} 2} \int_0^{+\infty} |h(\si)| d\si  \tag {2.56}
\end{align*}
and
\begin{align*}III_2^{(2)}  \le  & C  (bs)^2 \int_{a-b}^{a+b} |h(\si)| (\si
s)^m |\la'(\si s)| d\si \int_{\si s}^{(a+b)s} \f  {\la(t)} {\la^2(\si s)}dt \qquad (  \text{by (2.1)})\\
=& C (bs)^2  \int_{a-b}^{a+b} |h(\si)| (\si
s)^m \f {|\la'(\si s)|} {\la(\si s)} d\si \int_{\si s}^{(a+b)s} \f  {\la(t)} {\la(\si s)}dt\\
\le & C (bs)^3 \int_{a-b}^{a+b} |h(\si)| (\si s)^{ \f {3m} 2}  d\si
\qquad\qquad \qquad \qquad (  \text{by (2.1) and (2.4)})
\\ \le & C (bs)^3  (as)^{ \f {3 m} 2} \int_0^{+\infty} |h(\si)| d\si.  \tag {2.57}
\end{align*}
If $ (a-b)s \le 1 \le (a+b)s,$ then one has by (2.1)
\begin{align*} III_2^{(1)}
\le & C bs \int_{a-b}^{s^{-1} }|h(\si)| (\si s)^m |\la'(\si s)|
d\si \int_{(a-b)s}^{1} \la(t) dt \int_{(a-b)s}^t \f
{dy}{\la^2(y)}\\& + C bs \int_{s^{-1} }^{a+b }|h(\si)| (\si s)^m
|\la'(\si s)| d\si \int_{(a-b)s}^{\si s} \la(t) dt \int_{(a-b)s}^t
\f {dy}{\la^2(y)}\\
\le & C bs \int_{a-b}^{a+b }|h(\si)| (\si s)^m |\la'(\si s)| d\si+
C(bs)^3\int_{s^{-1} }^{a+b }|h(\si)| (\si s)^m
\f {|\la'(\si s)|} {\la(\si s)} d\si\\
\le & C \Big(bs (as)^{ \f m 2} + (bs)^3 (as)^{ \f {3m} 2} \Big)
\int_0^{+\infty} |h(\si)| d\si  \tag {2.58}
\end{align*}
and
\begin{align*} III_2^{(2)}
\le & C bs\int_{a-b}^{s^{-1}}|h(\si)| (\si s)^m |\la'(\si s)| d\si
\int_{\si s}^{(a+b)s} \la(t) dt \int_{(a-b)s}^{1} \f
{dy}{\la^2(y)}\\
& \  + C (bs)^2 \int_{s^{-1}}^{a+b}|h(\si)| (\si s)^m  \f{|\la'(\si
s)|} {\la(\si s)} d\si \int_{\si s}^{(a+b)s} \f {\la(t)} {\la(\si
s)} dt \\
\le & C \f 1 {\la^2(1)} (bs)^2\int_{a-b}^{s^{-1}}|h(\si)| \ (\si
s)^m \ |\la'(\si s)| \ \la(\si s) d\si +C (bs)^3
\int_{s^{-1}}^{a+b}|h(\si)| (\si s)^{ \f {3m}  2} d\si \\
\le & C \Big( (bs)^2 (as)^m +(bs)^3 (as)^{ \f {3m} 2}  \Big) \int_0^{+\infty}
|h(\si)| d\si. \tag {2.59}
\end{align*}
If $(a+b)s \le 1,$ then
\begin{align*} III_2^{(1)} \le & C bs \int_{a-b}^{s^{-1} }|h(\si)| (\si
s)^m |\la'(\si s)| d\si \int_{(a-b)s}^{1} \la(t) dt \int_{(a-b)s}^t
\f {dy}{\la^2(y)} \\
\le & C bs \int_{a-b}^{s^{-1} }|h(\si)| (\si s)^m |\la'(\si s)|
d\si \int_{(a-b)s}^{1} \f t {\la(t)} dt \qquad  (  \text{by (2.1)})
\\ \le & C bs (as)^{ \f m 2} \int_0^{+\infty} |h(\si)| d\si\tag {2.60}
\end{align*}
and
\begin{align*}  III_2^{(2)}
\le & C bs\int_{a-b}^{a+b}|h(\si)| (\si s)^m |\la'(\si s)| d\si
\int_{\si s}^{(a+b)s} \la(t) dt \int_{(a-b)s}^{1}
\f {dy}{\la^2(y)} \\
\le & C \la^{-2}(1)  (bs)^2 \int_{a-b}^{a+b}|h(\si)| (\si s)^m
|\la'(\si s)| \la(\si s) d\si \qquad (  \text{by (2.1)})
\\ \le & C(bs)^2 (as)^m \int_0^{+\infty} |h(\si)| d\si. \tag {2.61}
\end{align*}

Combining (2.55)-(2.61) yields
$$III_2 \le C\Big(bs (as)^{ \f m 2}
+ (bs)^2 (as)^m + (bs)^3 (as)^{ \f {3 m} 2} \Big)\int_0^{+\infty} |h(\si)| d\si.
\eqno{( 2.62)} $$

\vskip 0.3 true cm

{\bf Step 3. Estimate of $III_3$}

\vskip 0.2 true cm

A direct computation derives that from (2.3)
\begin{align*}  III_3
 \le & C
 \Big( bs \la((a-b)s) \Big) \Big(bs \big( (a+b)s \big)^{m-1}\la((a+b)s) \Big)
  \Big(\f{(a-b)s}{\la^2((a-b)s)} \Big)
\int_0^{+\infty} |h(\si)| d\si \\
\le & C (bs)^2 (as)^m \int_0^{+\infty} |h(\si)| d\si. \tag {2.63}
\end{align*}

\vskip 0.3 true cm
{\bf Step 4. Estimate of $III_4$}

\vskip 0.2 true cm

We have
\begin{align*} III_4
 \le & C bs \big((a+b)s \big)^{m-1}\la((a+b)s) \int_{a-b}^{a+b}
 |h(\si)| d\si\int_{(a-b)s}^{\si s} \la(t) dt \int_{(a-b)s}^{t}  \f {dy}{\la^2(y)} \\
 & \
 + C bs \big((a+b)s \big)^{m-1}\la((a+b)s) \int_{a-b}^{a+b}
 |h(\si)| d\si\int_{\si s}^{(a+b)s} \la(t) dt\int_{(a-b)s}^{\si s}  \f
 {dy}{\la^2(y)}\\
\equiv& III_4^{(1)} + III_4^{(2)}.
\end{align*}

It follows from Lemma 2.2 and (2.1) that
\begin{align*} III_4^{(1)}
 \le &  C bs \big((a+b)s \big)^{m-1}\la((a+b)s)
\int_{a-b}^{a+b}
 |h(\si)| d\si\int_{(a-b)s}^{\si s} \f{bs} {\la(t)} dt \\
 \le & C (bs)^2 \big((a+b)s \big)^{m-1} \int_{a-b}^{a+b}
 |h(\si)| (\si s)\f {\la((a+b)s)}  {\la(\si s)} d\si\\
 \le & C(bs)^2 (as)^{m} \int_0^{+\infty} |h(\si)| d\si
\end{align*}
and
\begin{align*} III_4^{(2)}
\le & C bs \big((a+b)s \big)^{m-1}\la((a+b)s)
\int_{a-b}^{a+b}
|h(\si)| d\si \Big(\int_{\si s}^{(a+b)s} \la(t)dt\ \f {\si s} {\la^2(\si s)} \Big)\\
\le & C bs\big((a+b)s \big)^{m-1}
\int_{a-b}^{a+b}
|h(\si)| (\si s) \f { \la((a+b)s)}  {\la(\si s)}  d\si\int_{\si s}^{(a+b)s} \f {\la(t)} {\la(\si s)}
dt \\
\le & C (bs)^2 (as)^{m}  \int_0^{+\infty} |h(\si)| d\si.
\end{align*}
Hence
$$III_4 \le C (bs)^2 (as)^m \int_0^{+\infty} |h(\si)| d\si. \eqno{(2.64)}$$

\vskip 0.3 true cm

{\bf Step 5. Estimate of $III_5$}

\vskip 0.2 true cm

One has
\begin{align*}  III_5
=& \big((a+b)s \big)^m \la((a+b)s) \int_{a-b}^{a+b} |h(\si)| d\si
\int_{\si s}^{(a+b)s} \la(t) dt \int_{\si s}^t \f {dy} {\la^2(y)}\\
\le & \big((a+b)s \big)^m \la((a+b)s) \int_{a-b}^{a+b} |h(\si)| d\si
\int_{\si s}^{(a+b)s} \f{2bs} {\la(t)}  dt \\
\le & C (bs)^2 (as)^m \int_0^{+\infty} |h(\si)| d\si. \tag {2.65}
\end{align*}

Substituting (2.54) and (2.62)-(2.65) into (2.53) derives the conclusion
in Part 3. Therefore, collecting all the results in Part 1 -Part 3, we have
$$\int_{0}^{+\infty} \Big| \int_0^{+\infty} \hat{K}(t,
\sigma, \xi) h( \si) d\si \Big| dt \le C P_2(as,bs)\int_0^{+\infty} |h(\si)|
d\si, \eqno{(2.66)}$$ where $P_2(a,b)=b a^{\f m 2} + b^2 a^m +b^3 a^{\f
{3m} 2}.$
Hence Lemma 2.7 is proved.\hfill $\square$

Based on Lemma 2.5-Lemma 2.7, we will show that the function $\hat K(t,\si,\xi)$ defined in (2.16)
satisfies (2.20)-(2.22) for all $|\al|\le q$ with some integer $q>n$.

\vskip 0.3 true cm

{\bf Lemma 2.8.} {\it  Let $\hat{K}(t, \sigma, \xi)$ be defined in (2.16), then (2.20)-(2.22) hold.}

\vskip 0.1 true cm

{\bf Proof.}  We now prove (2.20) by induction method. Note that (2.20) holds for $|\al|=0$
by Lemma 2.5. Suppose that for $|\al|\le j<q$
$$\int_0^{+\infty} \Big|\partial_\xi^\alpha \Big(\int_0^{+\infty} \hat{K}(t, \sigma,
\xi) h(\sigma) d\sigma \Big) \Big| dt \le C_q |\xi|^{-|\alpha|}\int
|h(\sigma)| d\sigma\quad\text{for any $h \in L^1(\R^+)$ };\eqno {(A_j)} $$
\begin{align*}
&\int_{C \Delta(a, 2^{q+1}b)} \Big|\partial_\xi^\alpha
\int_0^{+\infty} \hat{K}(t, \sigma, \xi) h(\sigma) d\sigma \Big|
dt \le C_q D_1(a,b, R)|\xi|^{-|\alpha|}\int |h(\sigma)| d\sigma\\
&\qquad\qquad\qquad\qquad\qquad\qquad\qquad
\text{for any $h\in L^1(\R^+)$ with $supp h\subset \Delta(a,b)$};
\tag {$B_j$}
\end{align*}
and
\begin{align*}
&\int_0^{+\infty} \Big|\partial_\xi^\alpha \int_0^{+\infty}
\hat{K}(t, \sigma, \xi) h(\sigma) d\sigma \Big| dt \le C_q
D_2(a,b, R)|\xi|^{-|\alpha|}\int |h(\sigma)| d\sigma\\
&\qquad\qquad \qquad  \text{for any $h \in L^1(\R^+)$ with $supp h\subset \Delta(a,b)$ and
$\int h(\si) d\si =0$}. \tag {$C_j$}
\end{align*}
At first, we show $(A_{j+1})$ holds. Denote
$U^{(j+1)}=\partial_\xi^{j+1}\hat{u}(t, \xi)$. Taking
$\p_{\xi}^{j+1}$ on the two hand sides of the equation
$\p_t^2\hat u-t^m|\xi|^2\hat u=t^m  h(t)$ yields
$$
 \p_t^2 U^{(j+1)}-t^m |\xi|^2 U^{(j+1)}=t^m \Big( C_1
\p_\xi(|\xi|^2)U^{(j)} +C_2 \p_\xi^2(|\xi|^2)U^{(j-1)}\Big),\eqno{(2.67)}
$$
where $C_1$ and $C_2$ are  some constants.
It follows from (2.67) and $(A_0)$ that
$$
\int_0^{+\infty} |\xi|^2 \ |U^{(j +1)}(t,\xi)| dt \le
C |\xi|^{-1}\int_0^{+\infty} |\xi|^2 \ |U^{(j)}(t,\xi)| \ dt +C
|\xi|^{-2}\int_0^{+\infty} |\xi|^2 \ |U^{(j-1 )}(t,\xi)| \ dt.\eqno{(2.68)}
$$
This, together with the induction assumption $(A_{j})$ and the definition of $\hat{K}(t, \sigma, \xi)$,
yields $(A_{j+1})$.

In addition, $(B_{j+1})$ or $(C_{j+1})$ can be directly obtained
by using (2.68) repeatedly until $j=1$ and combining Lemma 2.6 or Lemma 2.7.
Therefore, the proof of Lemma 2.8 is completed. \hfill $\square$

Next we show that functions $$D_1(a,b,r)=P_1(ar^{\f{2}{m+2}}, br^{\f{2}{m+2}}), \ \ D_2(a,b,r)=P_2(ar^{\f{2}{m+2}},
br^{\f{2}{m+2}})$$ satisfy the estimate (2.23) for suitable $a_k, b_k$ and $\t b_k$, where $P_1(a,b)$ and $P_2(a,b)$ are defined in (2.29)
and (2.33) respectively.

\vskip 0.3 true cm

{\bf Lemma 2.9.} {\it For $D_1(a,b,r)$ and $D_2(a,b,r)$ defined above,
the estimate (2.23) in Lemma 2.4 holds for $a_k, b_k$ and $\t b_k$
with $b_k\sim a_k^{m/2} \tilde{b}_k$.}

\vskip .1cm

{\bf Proof.} For any $f\in L^1(\Bbb R^+\times\Bbb R^n)$, one has the standard Calder\'on-Zygmund decomposition
(see Theorem 4 in Chapter 1 of [21]):
$$f=g+\ds\sum r_k$$
such that
\begin{equation*}
g(t,x)= \left\{
\begin{aligned}
 & f(t,x), \qquad \qquad \qquad \quad \text{in}  \ \R^{n+1}_+\setminus \cup_k Q_k \\
& \f 1 {|Q_k|}\int_{Q_k} f(\theta,y) d\theta dy, \qquad  \text{in} \ \ Q_k
\end{aligned}
\right.
\end{equation*}
and $$r_k(t,x) =\Big(f(t,x) -g(t,x)\Big)\chi_{Q_k},$$ where the cube
 $Q_k=\Delta_k(a_k, \tilde{b}_k) \times
I_k(x_k, b_k)$ with $\ds\f 1 C \le \ds\f {b_k}{a_k^{m/2} \tilde{b}_k} \le
C$ for some positive constant $C>1$ independent of $k$.

Set $$D_1(a_k,\t b_k,r)\equiv P_1(a_k
r^{2/(m+2)}, \t b_k r^{2/(m+2)} ) = exp\big(-C a_k^{\f m 2}\t b_k r\big),$$ and
$$
D_2(a_k,\t b_k,r)\equiv a_k^{\f m 2}\t b_k r+a_k^m\t b_k^2r^2+a_k^{\f {3m}{2}}\t b_k^3 r^3,$$
then \begin{align*} &\ds\sum_{2^jb_k\ge 1}D_1(a_k,\t b_k,2^j)+\ds\sum_{2^jb_k\le 1}D_2(a_k,\t b_k,2^j)\\
\le &\ds\sum_{2^jb_k\ge 1}exp (-C2^jb_k)+C\ds\sum_{2^jb_k\le 1}\bigg(2^jb_k+ (2^jb_k)^2+(2^jb_k)^3\biggr) \\
& \le C \end{align*}
and Lemma 2.9 is proved. \hfill $\square$
\vskip .3cm

Based on Lemma 2.4-Lemma 2.9, we now prove

\vskip .1cm

{\bf Theorem 2.10.} {\it Consider the problem
\begin{equation*}
\left\{
\begin{aligned} &\p_t^2w+t^{m}\Delta w=t^m g(t,x),\quad (t,x)\in [0,+\infty) \times\Bbb R^n,\\
&w(0,x)=0,
\end{aligned}
\right.\tag{2.69}
\end{equation*}
where $m \in \N$, $g \in L^p(\R^{n+1}_+)$ $(1<p<\infty)$, and
$supp g\subset\{(t,x): 0\le t\le M_0\}$ with $M_0>0$ being some fixed constant, then (2.69) has a unique solution
$w\in W^{2,p}([0, T]\times \Bbb R^{n})$ for any $T>0$, moreover $w$ satisfies the following estimate
\begin{align*}
&\| \p_t^2 w\|_{L^p(G_T)} + \|\Delta w\|_{L^p(G_T)}
+ \sum_{j=1}^n\|\p_{tx_j}^2 w\|_{L^p(G_T)} \\
& \ \ +\| \p_t w\|_{L^p(G_T)} +\sum_{j=1}^n\|\p_{j} w\|_{L^p(G_T)}
+\|w\|_{L^p(G_T)}\\
& \le C_{p,T}
\|g\|_{L^p(\R^{n+1}_+)}, \tag{2.70}
\end{align*} where $G_T=[0, T] \times \R^n$, and $C_{p,T}$ is a generic positive constant depending on $p$
and $T$. }

\vskip .1cm
{\bf Proof.}  By (2.16), we know that the solution $u$ to (2.69) satisfies
$$
\Delta w=  \mathcal F_\xi^{-1} \big(|\xi|^2 \hat{u}(t,\xi)
\big)=\mathcal F_\xi^{-1} \Big( \int_0^{+\infty} \hat{K}(t, \sigma,\xi)
\hat{g}(\sigma,\xi) d\sigma \Big),
$$
where $\hat{K}(t, \sigma, \xi)=|\xi|^2\sigma^m
\hat{T}(t,\sigma,\xi)$.
By Lemma 2.5-Lemma 2.9, one knows that  $\hat{K}(t, \sigma, \xi)$ satisfies all the requirements in Lemma 2.4. Hence we have from Lemma 2.4
$$
\|\Delta w\|_{L^p(\Bbb R_+^{n+1})}\le C_{p}\|g\|_{L^p(\Bbb R_+^{n+1})}.\eqno{(2.71)}$$
This, together with the equation (2.69), yields
$$\|\p_t^2w\|_{L^p(G_T)}\le C_p T^m \|g\|_{L^p(G_T)}.\eqno{(2.72)}$$
By (2.71)-(2.72) and the interpolation theory, we can obtain
$$\sum_{j=1}^n\|\p_{tx_j}^2 w\|_{L^p(G_T)} \le C_{p,T}
\|g\|_{L^p(\R^{n+1}_+)}.\eqno{(2.73)}$$
In addition, by (2.15) we have
$$\hat w(t,\xi)=\la(ts)\int_0^t\f{dy}{\la^2(ys)}\int_y^{\infty}\la(\si s)\si^m\hat g(\si, \xi)d\si,$$
which derives
$$\p_t\hat w(0,\xi)=\int_0^{\infty}\la(\si s)\si^m\hat g(\si, \xi)d\si=
\int_0^{M_0}\la(\si s)\si^m\hat g(\si, \xi)d\si$$
and further
$$\p_tw(0,x)=\int_0^{M_0}\mathcal F_\xi^{-1}(\la(\si s))\ast (\si^mg(\si, \cdot))d\si.\eqno{(2.74)}$$

It is noted that
$$\p_{\xi}^\al(\la(-\si s))=|\xi|^{-|\al|} \Big( P_1(\si s,\o) \la(-\si s)
+P_2(\si s,\o) \la'(-\si s) \Big),\eqno{(2.75)}$$ where $P_j(\theta,\o) (j=1,2)$ are
the polynomials of $\theta$,  and smoothly depend on the variable $\o\in S^{n-1}$.
From (2.75), we easily derive that for some integer $q>n$
$$\sup \limits_{t \in [-T, 0]} r^{2 |\alpha|-n} \int_{\f r 2 \le |\xi| \le 2 r}
|\p_{\xi}^\al(\la(-ts))|^2 d\xi \le C_q
 \qquad \text{ for all $r>0$ and $|\al|\le q$.} \eqno{( 2.76)}$$
Analogously, we have  for some integer $q>n$
$$\sup \limits_{t \in [-T, 0]} r^{2 |\alpha|-n} \int_{\f r 2 \le |\xi| \le 2 r}
|\p_{\xi}^\al(\la'(-ts))|^2 d\xi \le C_q
 \qquad \text{ for all $r>0$ and $|\al|\le q$.}  \eqno{( 2.77)}$$
Then it follows (2.76)-(2.77), H\"ormander's multiplier theorem (see Theorem 7.95 in [10]), Minkowski inequality  and (2.74) that
$$\|\p_tw(0,x)\|_{L^p(\Bbb R^n)}\le C_p M_0^{m+1-\f{1}{p}}\|g\|_{L^p(\Bbb R^{n+1}_+)}.\eqno{(2.78)}$$
This, together with (2.72), yields
\begin{align*}
\|\p_t w\|_{L^p(G_T)} \le & \|\p_t w(0, \cdot)\|_{L^p(G_T)}+ \Big( \int_0^T \Big\| \int_0^t \p_\tau^2 w(
\tau, \cdot) d\tau \Big\|_{L^p(\Bbb R^n)}^p\Big)^{1/p}    \\
\le & CT^{\f1p}(\|\p_t w (0, \cdot)\|_{L^p(\Bbb R^n)}+ \|\p_\tau^2 w(
\tau, x)\|_{L^p(G_T)}) \\ \le  & C_p T^{\f1p}(M_0^{m+1-\f{1}{p}}+T^{m})\|g\|_{L^p(\Bbb R^{n+1}_+)}.\tag{2.79}
\end{align*}
By (2.73) and $w(0,x)=0$, one has from  Minkowski inequality, H\"older's inequality and (2.79) that
$$\|w\|_{L^p(G_T)}\le T \|\p_t w\|_{L^p(G_T)} \le  C_p T^{1+\f{1}{p}}(M_0^{m+1-\f{1}{p}}+T^{m})\|g\|_{L^p(\Bbb R^{n+1}_+)}.\eqno{(2.80)}$$
In addition, it follows from $w(0,x)=0$ and (2.73)
that
$$\sum_{j=1}^n\|\p_{j} w\|_{L^p(G_T)}
\le C_{p,T}
\|g\|_{L^p(\R^{n+1}_+)}. \eqno{(2.81)} $$
Combining (2.71)-(2.73) and (2.79)-(2.81) yield  (2.70). \hfill $\square$

\vskip 0.4 true cm
\section{ Weighted $W^{2,p}$ estimates for the
 generalized Tricomi  equations}

In this section, for the later requirements of solving Theorem 1.1 in the degenerate elliptic region, we will establish
the weighted  $W^{2,p}$ estimates of the solution to the following problem
\begin{equation*}
\left\{
\begin{aligned} &\p_t^2w+t^{m}\Delta_x w=t^\nu g(t,x),\quad (t,x)\in [0,+\infty) \times\Bbb R^n,\\
&w(0,x)=0,
\end{aligned}
\right.\tag{3.1}
\end{equation*}
where $m \in \N$, $\nu \in \Bbb R$, and $0\le\nu\le m.$  Here we point out that showing the
weighted $W^{2,p}$ estimates for (3.1) also admits independent interests in the
study on the linear degenerate elliptic equations.
\vskip .3cm

{\bf Theorem 3.1.} {\it Let $g \in L^p(\R^{n+1}_+)$ with $1<p<\infty$ and $supp (g)\subset\{(t,x): 0\le t\le M_0\}$ for some $M_0>0$,
then for any fixed $T>0$, there exists a generic constant $C_{p,T}>0$ such that
\begin{align*}
&\| \p_t^2 w\|_{L^p(G_T)} + \|t^{m-\nu} \Delta w\|_{L^p(G_T)}
+ \sum_{j=1}^n\|t^{\f {m-\nu}{2}}\p_{tx_j}^2 w\|_{L^p(G_T)}\\
&\quad +\| \p_t w\|_{L^p(G_T)}
+\sum_{j=1}^n\|t^{\f {m-\nu}{2}}\p_{x_j} w\|_{L^p(G_T)}
+\|w\|_{L^p(G_T)}
\le C_{p,T}
\|g\|_{L^p(\R^{n+1}_+)},\tag{3.2}
\end{align*}
where $G_T=[0,T] \times \R^n$.
In addition, we have further that $w, \p_t w\in C([0,T], L^p(\R^n))$ and
$$ \|w(t, \cdot)\|_{L^p(\R^n)}+\| \p_t w(t, \cdot)\|_{L^p(\R^n)}\le C_{p,T} \|g\|_{L^p(\R^{n+1}_+)}.\eqno{(3.3)}$$
Moreover,  if  $1<p<Q_\nu\equiv 1+ n \big( \f {2l+1-\nu} 2 \big ),$   then $w\in C([0,T], L^{\bar p}(\R^n)),
\p_t w\in L^{\bar p}(G_T)$ with $\ds\f{1}{\bar p}=\f{1}{p}-\f{1}{Q_\nu}$,
and  $$\|w(t, \cdot)\|_{L^{\bar p}(\R^n)}+\|\p_t w\|_{L^{\bar p}(G_T)} \le C_{p,T}\|g\|_{L^p(\R^{n+1}_+)};\eqno{(3.4)}$$
if $p>Q_\nu$,
then $w\in C([0,T], L^\infty(\R^n)), \p_t w\in L^{\infty}(G_T)$, and
$$\|w(t, \cdot)\|_{L^{\infty}(\R^n)}+\|\p_t w\|_{L^{\infty}(G_T)} \le C_T \|g\|_{L^p(\R^{n+1}_+)}; \eqno{(3.5)}$$
and if $p=Q_\nu,$ then $w\in C([0,T], L^q(\R^n)), \p_tw\in L^{q}(G_T)$ for any $1<q < \infty$, and
$$\|w(t, \cdot)\|_{L^{q}(\R^n)}+\|\p_t w\|_{L^{q}(G_T)} \le C_{q,T} \|g\|_{L^p(\R^{n+1}_+)}, \eqno{(3.6)} $$
}

\vskip .2cm
{\bf Remark 3.1.} {\it When $\nu=0$  in (3.1), the weighted $W^{2,p}$ or $W^{1,p}$ estimates have been established under variable assumptions on the function $g(t,x)$ and by different methods in [11], [21], and [24].}

\vskip .2cm

{\bf Remark 3.2.} {\it By [17] and the references therein, $Q_\nu=1+(\f{m-\nu}{2}+1)n$ stands for the homogeneous dimension with respect to the operator
$\p_t^2+t^{m-\nu}\Delta$ for $t\ge 0$ and $x\in\Bbb R^n$, and then it follows from the ``Sobolev's imbedding theorem" and (3.2)
that (3.4)-(3.6) hold.}

\vskip .2cm

{\bf Proof of Theorem 3.1.}
Since $w$ is a solution to (3.1), then we have from (2.14)-(2.15) that
$$\hat{w}(t,\xi)=\int_0^{+\infty} \hat{T}(t, \sigma, \xi)  \sigma^\nu \hat{g}(\sigma,\xi) d\sigma, $$
where  $ \hat{T}(t, \sigma, \xi)=\ds\int_0^{\min(t, \sigma)} \f
{\la(ts) \la(\sigma s)} {\la^2(ys)} dy$ and  $s= |\xi|^{2/{m+2}}$. Thus,
$$w(t,x)=\int_0^{+\infty} T(t, \sigma,\cdot) \ast \big(\sigma^\nu g(\sigma,\cdot)\big) d\sigma, \eqno{( 3.7)}$$
where $T(t,\sigma,x)=\mathcal F^{-1}_\xi \hat{T}(t, \sigma, \xi)$. To
show Theorem 3.1, as in Theorem 2.10, at first, we require to prove
$$\|t^{m-\nu}\Delta
w\|_{L^p(G_T)} \le C \|g\|_{L^p}. \eqno{( 3.8)} $$
It is noted that
$$
t^{m-\nu}\Delta w=  \mathcal F_\xi^{-1} \big( t^{m-\nu} |\xi|^2 \hat{w}(t,\xi)
\big)=\mathcal F_\xi^{-1} \Big( \int_0^{+\infty} \hat{K}(t, \sigma,\xi)
\hat{g}(\sigma,\xi) d\sigma \Big),
$$
where $\hat{K}(t, \sigma, \xi)=|\xi|^2 t^{m-\nu} \sigma^\nu
\hat{T}(t,\sigma,\xi)$. For notational convenience, we set
$\hat{K}_1(t, \sigma, \xi)=|\xi|^2 \sigma^m
\hat{T}(t,\sigma,\xi)$
and $\hat{K}_2(t, \sigma, \xi)=|\xi|^2 t^{m}\hat{T}(t,\sigma,\xi)$.

With respect to $\hat{K}_1(t, \sigma, \xi)$, its some crucial properties have been established in Lemma 2.5-Lemma 2.9
of $\S 2$.
For $\hat{K}_2(t, \sigma, \xi)$, the authors in [8] have shown that it satisfies (2.17)-(2.20) and (2.21)-(2.22)
with the quantities
$$D_1(a,b,r)=P_3(a r^{2/(m+2)}, b r^{2/(m+2)} ), \ \  D_2(a,b,r)=P_4(a r^{2/(m+2)}, b r^{2/(m+2)}),$$
where
\begin{equation*}
P_3(a,b)=\left\{
\begin{aligned}
  & a \qquad \qquad \qquad \qquad\text{ if } \ a \le 3, \\
& exp(-Ca^{m/2} b) \qquad\text{ if } \ a>3,
\end{aligned}
\right.\tag{3.9}
\end{equation*}
and
\begin{equation*}
P_4(a,b)=\left\{
\begin{aligned}
  & a \qquad \qquad \qquad \qquad \qquad\text{ if } \ a \le 3, \\
&  a^{m/2}b +  a^m b^2 + a b^{2/m}\quad \ \text{ if } \ a>3.
\end{aligned}
\right.\tag{3.10}
\end{equation*}

We now prove that $\hat{K}(t, \sigma, \xi)$ satisfies (2.17)-(2.20) and (2.21)-(2.22)
with suitable quantities $D_1(a,b,r)$ and $D_2(a,b,r)$. Notice that
\begin{equation*}
y^{m-\nu}\sigma^\nu \le \left\{
\begin{aligned}
 & y^m  \qquad \text{ if } \si \le y,\\
& \si^m \qquad \text{ if } y \le \si,
\end{aligned}
\right.
\end{equation*}
then one has $0\le \hat{K}(t, \sigma, \xi) \le \hat{K}_1(t, \sigma, \xi)
+\hat{K}_2(t, \sigma, \xi)$. From this, we have
$$\sup_{t, \xi} \int_0^{+\infty} |\hat{K}(t, \sigma, \xi)| d\sigma \le
\sup_{t, \xi} \int_0^{+\infty} |\hat{K}_1(t, \sigma, \xi)| d\sigma + \sup_{t,
\xi} \int_0^{+\infty} |\hat{K}_2(t, \sigma, \xi)| d\sigma \le C,\eqno{(3.11)}$$
$$\sup_{\si, \xi} \int_0^{+\infty} |\hat{K}(t, \sigma, \xi)| d\sigma \le
\sup_{\si, \xi} \int_0^{+\infty} |\hat{K}_1(t, \sigma, \xi)| dt + \sup_{\si,
\xi} \int_0^{+\infty} |\hat{K}_2(t, \sigma, \xi)| dt \le C,\eqno{(3.12)}$$
\begin{align*}   \int_0^{+\infty} \Big| \int_0^{+\infty} \hat{K}(t,
\sigma, \xi) h(\sigma) d\sigma \Big| dt  \le & \int_0^{+\infty}dt \int_0^{+\infty} \Big(\hat{K}_1(t, \sigma,
\xi)+ \hat{K}_2(t, \sigma, \xi)\Big)|h(\si)| d\sigma\\
\le & C \int_0^{+\infty} |h(\si)| d\si\tag{3.13}
\end{align*}
and
\begin{align*}  &  \int_{C \Delta(a, 2^{q+1}b)} \Big| \int_0^{+\infty}
\hat{K}(t, \sigma, \xi) h(\sigma) d\sigma \Big| dt \\ \le & \int_{C \Delta(a, 2^{q+1}b)}dt \int_0^{+\infty}
\Big(\hat{K}_1(t, \sigma, \xi)+\hat{K}_2(t, \sigma, \xi) \Big)
|h(\sigma)| d\sigma\\
\le & C  \Big( P_1(as,bs)+ P_3(as,bs) \Big) \int_0^{+\infty} |h(\si)| d\si,\tag{3.14}
\end{align*}
which means that $\hat{K} $ satisfies the estimates (2.18)-(2.19),
and (2.20)-(2.21) with $\alpha=0$.

Next we verify the
estimate (2.22) of $\hat{K}$ when $|\al|=0$. This procedure will be divided into
the following three parts. From now on, we assume that the function $h \in C^\infty_0(\R^+)$
with $\int h(\si)d\si =0$ is supported in
$\Delta(a,b)$.

{\bf Part 1.  $\int_{a+b}^{+\infty} \Big| \int_0^{+\infty} \hat{K}(t,
\sigma, \xi) h( \si) d\si \Big| dt \le C P_5(as, bs) \int_0^{+\infty} |h(\si)|
d\si$ holds for suitable function $P_5(a, b)$}

\vskip 0.3 true cm

It follows from a direct computation that
\begin{align*}  & \int_{a+b}^{+\infty} \Big| \int_0^{+\infty} \hat{K}(t,
\sigma, \xi) h( \si) d\si \Big| dt \\
=& \int_{a+b}^{+\infty} \Big| \int_{a-b}^{a+b} \big( \hat{K}(t,
\sigma, \xi)- \hat{K}(t, a+b, \xi)  \big) h(\si) d\si   \Big| dt
\qquad  (\text{by}\int h(\si)d\si =0)
\\
=& \int_{a+b}^{+\infty} \Big| \int_{a-b}^{a+b} h(\si) |\xi|^2
t^{m-\nu}\Big( \si^\nu \int_0^{\min(t,\sigma)} \f { \la(ts) \la(\sigma s)}
{ \lambda^2(ys)} dy \\
& \qquad \qquad \qquad \qquad \qquad \qquad-  (a+b)^\nu \int_0^{\min(t,a+b)} \f { \la(ts)
\la((a+b) s)} { \lambda^2(ys)} dy
\Big)d\si \Big| dt \\
\le & \int_{(a+b)s}^{+\infty} \la(t) dt \int_{a-b}^{a+b} |h(\si)| \
t^m  \Big( \la(\si s) -\la((a+b)s) \Big) d\si \int_0^{\si s} \f
{dy}{\la^2(y)} \\
& \   + \int_{(a+b)s}^{+\infty} \la(t) dt \int_{a-b}^{a+b} |h(\si)|
\ t^{m-\nu} \Big(\big((a+b)s\big)^\nu - \big(\si s\big)^\nu \Big)
\la((a+b)s)
 d\si \int_0^{\si s} \f {dy}{\la^2(y)} \\
 & \ +\int_{(a+b)s}^{+\infty} \la(t) dt \int_{a-b}^{a+b} |h(\si)| \
 t^m \la((a+b)s) d\si \int_{\si s}^{(a+b)s} \f {dy}
 {\la^2(y)} \\
\equiv & I_1 +I_2 +I_3.\tag{3.15}
\end{align*}
Next we treat all the $I_i$ ($1\le i\le 3$) in the distinct two cases.

\vskip 0.2 true cm

{\bf Case 1. $a s \ge 1$}
\vskip 0.2 true cm

In this case, we arrive at
\begin{align*}  I_1
\le & |\la'((a+b)s)| \int_{a-b}^{a+b} |h(\si)| d\si \int_0^{(a+b) s}
\f {\la(\si s) -\la((a+b)s) }{\la^2(y)} dy \qquad (  \text{by } t^m \la(t)= \la''(t)) \\
\le & C bs |\la'((a+b)s)| \int_{a-b}^{a+b} |h(\si)| d\si
\int_0^{(a+b) s} \f {-\la'(y)} {\la^2(y)} dy \qquad  (  \text{  by (2.1) }) \\
=& C bs \f {|\la'((a+b)s)|} {\la((a+b)s)} \Big(1-\la((a+b)s) \Big)
\int_{a-b}^{a+b} |h(\si)| d\si \\
\le & C bs (as)^{\f m 2} \int_0^{+\infty} |h(\si)| d\si  \tag {3.16}
\end{align*}
and
\begin{align*}  I_2 \le & C bs \big((a+b)s \big)^{\nu-1}\la((a+b)s)\int_{a-b}^{a+b}
|h(\si)| \   \f {\si s}{\la^2(\si s)} d\si
\int_{(a+b)s}^{+\infty} t^{m-\nu} \la(t) dt \\
\le & C bs \int_{a-b}^{a+b} |h(\si)| \  \f {\la((a+b)s)}{\la^2(\si
s)} d\si
\int_{(a+b)s}^{+\infty} t^m \la(t) dt \\
= & C bs \int_{a-b}^{a+b} |h(\si)| \  \f {|\la'((a+b)s)| \ \la((a+b)s)}{\la^2(\si
s)} d\si \\
\le & C bs \int_{a-b}^{a+b} |h(\si)| \  \f
{|\la'((a+b)s)|}{\la((a+b) s)}  \f {\la^2((a+b)s)}{\la^2(\si s)}
d\si
\\ \le & C bs (as)^{ \f m 2} \int_0^{+\infty} |h(\si)| d\si \tag {3.17}
\end{align*}
and
\begin{align*}  I_3
\le  & C \f {bs} {\la((a+b)s)} \int_{(a+b)s}^{+\infty} t^m\la(t) dt
\int_0^{+\infty} |h(\si)| d\si \qquad  (  \text{by (2.1)})\\
=& C bs \f{|\la'((a+b)s)|} {\la((a+b)s)} \int_0^{+\infty} |h(\si)| d\si \\
\le & C bs (as)^{ \f m 2} \int_0^{+\infty} |h(\si)| d\si. \qquad
\qquad (  \text{by (2.2)}
 ) \tag {3.18}
\end{align*}

{\bf Case 2. $as \le 1$}

\vskip 0.2 true cm

Denote by $H=\ds\int_{a-b}^{a+b} |h(\si)| \la (\si s) d\si \int_0^{\si s}
\f{dy}
{\la^2(y)}  \int_{y}^{+\infty} t^m \la(t) dt$. Then a direct computation yields
\begin{align*}
H= &\int_{a-b}^{a+b} |h(\si)| \la (\si s) d\si \int_0^{\si s} \f
{-\la'(y)}{\la^2(y)} dy  \qquad ( \text { since } t^m
\la(t)=\la''(t) ) \\
=& \int_{a-b}^{a+b} |h(\si)| \Big( \la(0)-\la( \si s) \Big) d\si \\
\le & (a+b)s |\la'(0)| \int_0^{+\infty} |h(\si)| d\si   \qquad \qquad \qquad \qquad \qquad (
\text{  by (2.1) }) \\
\le & C as\int_0^{+\infty} |h(\si)| d\si.
\end{align*}
Thus,

$$
I_1\le H\le C as\int_0^{+\infty} |h(\si)| d\si,\qquad
I_2\le  H \le C as\int_0^{+\infty} |h(\si)| d\si.\eqno{( 3.19)}
$$
In addition,
\begin{align*} I_3 \le & \int_{a-b}^{a+b}
|h(\si)| \ \la((a+b)s) d\si \int_{0}^{(a+b)s} \f {dy}
 {\la^2(y)}  \int_{y}^{+\infty} t^m \la(t) dt
 \\
= & \int_{a-b}^{a+b} |h(\si)| \la ((a+b)s ) d\si \int_{0}^{(a+b)s}
\f {-\la'(y)}{\la^2(y)} dy\\
\le & C as\int_0^{+\infty} |h(\si)| d\si. \tag {3.20}
\end{align*}
Substituting (3.16)-(3.20) into (3.15) yields
$$ \int_{a+b}^{+\infty} \Big| \int_0^{+\infty} \hat{K}(t,
\sigma, \xi) h( \si) d\si \Big| dt \le C P_5(as, bs) \int_0^{+\infty} |h(\si)|
d\si, \eqno{( 3.21)}
$$
where
\begin{equation*}
P_5(a, b)=\left\{
\begin{aligned}
  & a \qquad \qquad \text { if } a \le 1 \\
& b a^{m/2} \qquad \text { if } a \ge 1
\end{aligned}
\right.\tag{3.22}
\end{equation*}

{\bf Part 2. $\int_0^{a-b} \Big| \int_0^{+\infty} \hat{K}(t, \sigma,
\xi) h(\si) d\si \Big| dt\le C bs (as)^{m/2} \int_0^{+\infty} |h(\si)| d\si$ holds}

It follows from a direct computation that
\begin{align*}  & \int_0^{a-b} \Big| \int_0^{+\infty} \hat{K}(t, \sigma,
\xi) h(\si) d\si \Big| dt \\
=& \int_0^{a-b} t^{m-\nu}\la(ts) dt\Big| \int_{a-b}^{a+b} h(\si) |\xi|^2
\
 \Big( \si^\nu \int_0^{t} \f {  \la(\sigma s)} { \lambda^2(ys)}
dy - (a+b)^\nu \int_0^{t} \f {  \la((a+b) s)} { \lambda^2(ys)} dy
\Big) d\si \Big| \\
=& \int_0^{(a-b)s} t^{m-\nu} \la(t) dt\Big| \int_{a-b}^{a+b} h(\si) \
\Big( (\si s)^\nu  \la(\sigma s) \int_0^{t} \f {dy } { \lambda^2(y)} -
((a+b)s)^\nu \la((a+b) s)\int_0^{t} \f {dy } { \lambda^2(y)}  \Big)
d\si \Big|
\\
\le &   \int_{a-b}^{a+b} |h(\si)| \ (\si s)^m \Big( \la(\si s)
-\la((a+b)s) \Big)d\si \int_0^{(a-b)s} \la(t) dt \int_0^t \f {dy} {
\lambda^2(y)} \\
& \ +  \la((a+b)s) \int_{a-b}^{a+b} |h(\si)|  \Big((\si s)^\nu
-((a+b)s)^\nu \Big)d\si  \int_0^{(a-b)s}t^{m-\nu} \la(t) dt \int_0^t \f
{dy} {
\lambda^2(y)} \\
\le & C bs\int_{a-b}^{a+b} |h(\si)|   (\si s)^m |\la'(\si s)| d\si
\int_0^{(a-b)s} \la(t) dt \int_0^t \f {dy} { \lambda^2(y)}\\
& \ +C bs ((a+b)s)^{m-1} \la((a+b)s) \int_{a-b}^{a+b} |h(\si)|
d\si \int_0^{(a-b)s} \la(t) dt \int_0^t \f {dy} { \lambda^2(y)}
\qquad  (\text{by
(2.1)})\\
\le & C bs (as)^{m/2} \int_0^{+\infty} |h(\si)| d\si.\qquad\qquad\qquad\qquad\qquad\qquad\qquad\qquad
\qquad \qquad\quad \text{( as in (2.38) )} \tag {3.23}
\end{align*}

{\bf Part 3. $\ds\int_{a-b}^{a+b} \Big| \int_0^{+\infty} \hat{K}(t,
\sigma, \xi) h( \si) d\si \Big| dt  \le  C (bs)^2 (as)^m \int_0^{+\infty}
|h(\si)|d\si$ holds}

We have
\begin{align*}  & \int_{a-b}^{a+b} \Big| \int_0^{+\infty} \hat{K}(t,
\sigma, \xi) h( \si) d\si \Big| dt\\
=& \int_{a-b}^{a+b} \Big| \int_{a-b}^{a+b} h(\si) |\xi|^2 \ t^{m-\nu}
\Big( \si^\nu \int_0^{\min(t,\sigma)} \f { \la(ts) \la(\sigma s)} {
\lambda^2(ys)} dy -  (a+b)^\nu \int_0^{t} \f { \la(ts) \la((a+b) s)} {
\lambda^2(ys)} dy
\Big)d\si \Big| dt \\
\le &
\int_{a-b}^{a+b} |h(\si)| (\si s)^m \Big( \la(\si s)
-\la((a+b)s) \Big) d\si \int_{(a-b)s}^{\si s} \la(t)
dt\Big(\int_0^{(a-b)s}+ \int_{(a-b)s}^{t} \Big)\f
{dy}{\la^2(y)} \\
& \ +\int_{a-b}^{a+b} |h(\si)|  \Big( \la(\si s) -\la((a+b)s) \Big)
d\si \int_{\si s }^{(a+b)s} t^m\la(t) dt\Big(\int_0^{(a-b)s}+
\int_{(a-b)s}^{\si s} \Big)\f
{dy}{\la^2(y)} \\
& \   +  \int_{a-b}^{a+b} |h(\si)| \Big(\big((a+b)s\big)^\nu -
\big(\si s\big)^\nu \Big)
\la((a+b)s) d\si  \int_{(a-b)s}^{\si s} t^{m-\nu}\la(t) dt\Big(\int_0^{(a-b)s}+\int_{(a-b)s}^{t} \Big) \f {dy}{\la^2(y)} \\
& \   + \int_{a-b}^{a+b} |h(\si)| \Big(\big((a+b)s\big)^\nu - \big(\si
s\big)^\nu \Big)
\la((a+b)s) d\si  \int_{\si s}^{(a+b)s} t^{m-\nu}\la(t) dt \Big(\int_0^{(a-b)s}+\int_{(a-b)s}^{\si s} \Big) \f {dy}{\la^2(y)} \\
& \ +\la((a+b)s) \int_{a-b}^{a+b} |h(\si)|
d\si \int_{\si s}^{(a+b)s} t^m\la(t) dt\int_{\si s}^{t} \f {dy}
{\la^2(y)} \\
\equiv & III_1 + III_2 +V_1 +V_2 +V_3+V_3+V_4+V_5+V_6+V_7,\tag{3.24}
\end{align*}
where $III_1$ and $III_2$ have been defined in (2.53), whose estimates have been established in (2.54)
and (2.62) respectively. Next we treat each $V_i$ ($1\le i\le 7$). For the term $V_1$, we have
\begin{align*}  V_1
\le & \int_{a-b}^{a+b} |h(\si)|     d\si \int_{(a-b) s }^{(a+b)s}
t^m\la(t) dt
\int_0^{(a-b)s}\f {\la(\si s) -\la((a+b)s)}{\la^2(y)}dy \\
\le & C bs \int_{a-b}^{a+b} |h(\si)|     d\si \int_{(a-b) s
}^{(a+b)s} t^m\la(t) dt\int_0^{(a-b)s}\f {-\la'(y)}{\la^2(y)}dy\\
\le & C (bs)^2 \big((a+b)s \big)^m  \la((a-b) s)\int_{a-b}^{a+b} |h(\si)|    d\si
\Big( \f 1{\la((a-b)s)} -1\Big) \qquad (  \text{by (2.1)}) \\
\le & C(bs)^2 (as)^m   \int_0^{+\infty} |h(\si)|d\si. \tag {3.25}
\end{align*}
For the term $V_2$,
\begin{align*}  V_2=&\int_{a-b}^{a+b} |h(\si)|   d\si \int_{\si s }^{(a+b)s} t^m\la(t)
dt \int_{(a-b)s}^{\si s} \f
{\la(\si s) -\la((a+b)s)}{\la^2(y)}dy\\
\le & C (bs)^2 \big( (a+b)s\big)^m \int_{a-b}^{a+b} |h(\si)| \
\la(\si s) d\si \int_{(a-b)s}^{\si s} \f {-\la'(y)}{\la^2(y)} dy\\
=& C (bs)^2 \big( (a+b)s\big)^m \int_{a-b}^{a+b} |h(\si)| \ \la(\si
s) \Big( \f 1 {\la(\si s)} -\f 1 {\la((a-b)s)} \Big) d\si \\
\le & C (bs)^2 (as)^m \int_0^{+\infty} |h(\si)|d\si. \tag {3.26}
\end{align*}
For the term $V_3$,
\begin{align*}  V_3
\le & C bs\big( (a+b)s \big)^{\nu-1}\la((a+b)s) \int_{a-b}^{a+b}
|h(\si)| d\si \ \Big( bs\ (\si s)^{m-\nu} \la((a-b)s) \Big) \ \Big(
\f {(a-b)s}{\la^2((a-b)s)} \Big)\\
\le & C (bs)^2 \big((a+b)s \big)^m \f{\la((a+b)s)}{\la((a-b)s)} \int
|h(\si)| d\si \\
\le & C (bs)^2 (as)^m \int_0^{+\infty} |h(\si)| d\si. \tag {3.27}
\end{align*}
For the term $V_4$,
\begin{align*}  V_4
\le & C bs \big((a+b)s\big)^{m-1}\la((a+b)s) \int_{a-b}^{a+b}
|h(\si)| d\si \int_{(a-b)s}^{\si s} \la(t) dt \int_{((a-b)s)}^t \f {dy}{\la^2(y)}\\
\le & C bs \big((a+b)s\big)^{m-1} \la((a+b)s)\int_{a-b}^{a+b}
|h(\si)| d\si \int_{(a-b)s}^{\si s} \f {bs}{\la(t)} dt \\
\le & C (bs)^2 \big((a+b)s\big)^{m-1} \int_{a-b}^{a+b}
\si s|h(\si)| \f{\la((a+b)s)}{\la(\si s)}d\si \\
\le & C (bs)^2 (as)^m \int_0^{+\infty} |h(\si)| d\si. \tag {3.28}
\end{align*}
For the term $V_5$,
\begin{align*}  V_5
\le &  C (bs)^2 \big( (a+b)s\big)^m \int_{a-b}^{a+b} |h(\si)| \f
{\la((a+b)s) \la(\si s)} {\la^2((a-b)s)} d\si \qquad (  \text{by (2.1)}) \\
\le & C (bs)^2 (as)^m \int_0^{+\infty} |h(\si)| d\si. \tag {3.29}
\end{align*}
For the term $V_6$,
\begin{align*}  V_6
\le & C (bs)^2 \big( (a+b)s\big)^{m-1} \int_{a-b}^{a+b} |h(\si)|
\la((a+b)s) \f {\si s} {\la(\si s)}d\si \qquad (
\text{by (2.1)}) \\
\le & C (bs)^2 (as)^m \int_0^{+\infty}|h(\si)|d\si. \tag {3.30}
\end{align*}
Finally, for the term $V_7$,
\begin{align*}  V_7
 \le & C bs\la((a+b)s) \int_{a-b}^{a+b} |h(\si)| d\si \int_{\si s}^{(a+b)s} \f {t^m} {\la(t)}
 dt \\
 \le & C (bs)^2 \big( (a+b)s\big)^m \int |h(\si)|d\si \\
 \le & C (bs)^2 (as)^m \int_0^{+\infty} |h(\si)|d\si. \tag {3.31}
\end{align*}
Collecting (3.24)-(3.31) yields
$$
 \int_{a-b}^{a+b} \Big| \int_0^{+\infty} \hat{K}(t,
\sigma, \xi) h( \si) d\si \Big| dt  \le  C (bs)^2 (as)^m \int_0^{+\infty}
|h(\si)|d\si.$$
Thus, by Part 1-Part 3, we arrive at
$$
 \int_{0}^{+\infty} \Big| \int_0^{+\infty} \hat{K}(t,
\sigma, \xi) h( \si) d\si \Big| dt  \le  C P_6(as, bs) \int_0^{+\infty}
|h(\si)|d\si, \eqno{(3.32)}$$ where
$$P_6(a,b)=P_5(a,b)+b a^{m/2}+b^2 a^m.$$
Therefore, $\hat{K} $ satisfies the estimate (2.12) for $|\alpha|=0$.
Applying the similar arguments as done in Lemma 2.8, we can obtain
that $\hat{K} $ satisfies estimates (2.20)-(2.22) for any
$|\alpha| \ge 0$ by induction method.

In addition, it is easy to know that (2.23) holds for suitably chosen $a_k, b_k$ and $\t b_k$
with $b_k\sim a_k^{m/2} \tilde{b}_k$ when we set
$D_1(a,b,r)=P_1(a
r^{2/(m+2)}, b r^{2/(m+2)} ) + P_3(a r^{2/(m+2)}, b r^{2/(m+2)} )$
and
$D_2(a,b,r)=P_6(a
r^{2/(m+2)},$ $b r^{2/(m+2)}).$

Hence, applying Lemma 2.4
we could get
$$ \|t^{m-\nu} \Delta w\|_{L^p(\Bbb R^{n+1}_+)}
\le  C_p\|g\|_{L^p(\Bbb R^{n+1}_+)},\eqno{(3.33)} $$
and then as arguing in Theorem 2.10, we have
$$\| \p_t^2 w\|_{L^p(G_T)} \le  C_p T^{\nu} \|g\|_{L^p(\Bbb R^{n+1}_+)},  \qquad \sum_{j=1}^n\|t^{\f {m-\nu}{2}}\p_{tx_j}^2 w\|_{L^p(G_T)}\le  C_{p,T}\|g\|_{L^p}.\eqno{(3.34)}$$

Note that from (3.3)
$$\hat w(t,\xi)=\la(ts)\int_0^t\f{dy}{\la^2(ys)}\int_y^{\infty}\la(\si s)\si^\nu\hat g(\si, \xi)d\si,$$
which derives
$$\p_t\hat w(0,\xi)=\int_0^{\infty}\la(\si s)\si^\nu\hat g(\si, \xi)d\si=
\int_0^{M_0}\la(\si s)\si^\nu\hat g(\si, \xi)d\si.$$
Thus,
$$\p_tw(0,x)=\int_0^{M_0}\mathcal F_\xi^{-1}(\la(\si s))\ast (\si^\nu g(\si, \cdot))d\si,$$
and $\|\p_t w(0,x)\|_{L^p(\Bbb R^n)}\le C_{p}M_0^{\nu+1-\f{1}{p}}\|g\|_{L^p(\Bbb R^{n+1}_+)}$ holds as in the proof of Theorem 2.10. Together with
(3.34), this yields for $0\le t\le T$
$$\|\p_t w\|_{L^p(G_T)}\le C_pT^{\f1p}(M_0^{\nu+1-\f{1}{p}}+T^{\nu})\|g\|_{L^p(\Bbb R^{n+1}_+)}\eqno{(3.35)}$$
and
\begin{align*}
&\|\p_t w(t, \cdot)\|_{L^p(\Bbb R^n)}\le \|\p_t w(0, \cdot)\|_{L^p(\Bbb R^n)}+\int_0^t\|\p_{\tau}^2w(\tau,\cdot)\|_{L^p(\Bbb R^n)}d\tau\\
&\le C_{p}M_0^{\nu+1-\f{1}{p}}\|g\|_{L^p(\Bbb R^{n+1}_+)}+T^{1-\f1p}\|\p_{\tau}^2 w(\tau,x)\|_{L^p(G_T)}\\
&\le C_p (M_0^{\nu+1-\f{1}{p}}+T^{\nu+1-\f1p})\|g\|_{L^p(\Bbb R^{n+1}_+)}.\tag{3.36}
\end{align*}
From (3.35)-(3.36) and $w(0,x)=0$, one has
\begin{equation*}
\left\{
\begin{aligned}
&\|w(t, \cdot)\|_{L^p(\Bbb R^n)} \le C_p T(M_0^{\nu+1-\f{1}{p}}+T^{\nu+1-\f1p})\|g\|_{L^p(\Bbb R^{n+1}_+)},\\
&\|w\|_{L^p(G_T)}  \le C_pT^{1+\f1p}(M_0^{\nu+1-\f{1}{p}}+T^{\nu})\|g\|_{L^p(\Bbb R^{n+1}_+)}.
\end{aligned}
\right.\tag{3.37}
\end{equation*}
In addition, it follows from $w(0,x)=0$, (3.34) and Hardy's inequality
that
$$\sum_{j=1}^n\|t^{\f {m-\nu}{2}}\p_{j} w\|_{L^p(G_T)}
\le C_{p,T}
\|g\|_{L^p(\R^{n+1}_+)}. \eqno{(3.38)}$$
Combining (3.33)-(3.38) yields  (3.2)-(3.3). In addition, (3.4)-(3.6) can be derived
as indicated in Remark 3.2. \hfill $\square$

\vskip 0.4 true cm

\section{ Solvability  of (1.1) in the degenerate elliptic region $\{t\le 0\}$ }

In this section, we will establish the solvability and regularity of problem (1.1) in the degenerate elliptic region $\{ t \le 0\}$ by applying the weighted
$W^{2,p}$ estimates given in Theorem 3.1. Our main result is:

\vskip .2cm

{\bf Theorem 4.1.} {\it Under  the assumptions (1.2)-(1.3),
there exists a constant $T_0>0$ such that when $0 \le \mu \le 1$, or when $1 < \mu <p_0$ with $Q_0 \le \ds\f{p_0}{\mu-1}$, the following degenerate elliptic equation
\begin{equation*}
\left\{
\begin{aligned} &\p_t^2 u-t^{2l-1}\Delta u= f(t,x,u),\qquad (t,x)\in (-\infty, 0] \times\Bbb R^n,\\
&u(0,x)=\varphi(x) \in H^s(\R^n),
\end{aligned}
\right.\tag{4.1}
\end{equation*}
has a unique local solution $u$ satisfying
$$u(t,x)\in C\Big([-T_0, 0], L^{p_0}(\Bbb R^n)   \Big), \ \ \p_t u(t,x)\in C\Big([-T_0, 0], L^{p_1}(\Bbb R^n)\Big), \eqno{ (4.2)} $$
where $0 \le s < \ds\f n 2$, $p_0=\ds\f{2n}{n-2s}$, and $\ds\f{1}{p_1}=\f{1}{2}-\f{1}{n}\bigl(s-\f{2}{2l+1}\bigr)$.}

\vskip .1cm

{\bf Proof.}
 Without loss of generality, we assume that $\mu \ne 0.$
Let $\bar u(t,x)$ be a solution to the following linear problem
\begin{equation*}
\left\{
\begin{aligned} &\p_t^2\bar u-t^{2l-1}\Delta \bar u=0,\qquad (t,x)\in (-\infty, 0] \times\Bbb R^n,\\
&\bar u(0,x)=\varphi(x)  \in H^s(\Bbb R^n),  \qquad  x \in \Bbb R^n.
\end{aligned}
\right.
\end{equation*}
Then from (2.11), we have
$$
\bar u(t,x)\in C((-\infty, 0], H^s(\Bbb R^n)) \cap C^1((-\infty, 0], H^{s-\f{2}{2l+1}}(\Bbb R^n)), $$
which, together with Sobolev's embedding theorem, yields
$$ \bar u(t,x)\in C((-\infty, 0], L^{p_0}(\Bbb R^n)) \cap C^1((-\infty, 0], L^{p_1}(\Bbb R^n)).\eqno{( 4.3 )}$$

Choosing a smooth function $\chi(t)$ such that $\chi(t)\equiv 1$ for $t\ge -1$ and $\chi(t)\equiv 0$ for $t\le -2$.
For suitably small fixed constant $T_0>0$, we set $\chi_{T_0}(t)=\chi(\ds\f{t}{T_0})$.
Let $\bar u_{T_0}$ be a solution to the linear equation
\begin{equation*}
\left\{
\begin{aligned} &\p_t^2\bar u_{T_0}-t^{2l-1}\Delta \bar u_{T_0}=\chi_{T_0}(t)f(t,x,\bar u),\qquad (t,x)\in (-\infty, 0] \times\Bbb R^n,\\
&\bar u_{T_0}(0,x)=0.
\end{aligned}
\right.\tag{4.4}
\end{equation*}
Due to the compact support property of $f(t,x,\bar u)$ on the variable $x$ and the regularity of $\bar u(t,x)$
in (4.3), one has from (1.2) that
$\chi_{T_0}(t)f(t,x,\bar u) \in L^p(\Bbb R^{n+1}_-)$, here $1< p \le \ds\f{p_0}{\mu}$
and $\Bbb R^{n+1}_-=\{(t,x): t\le 0, x\in\Bbb R^n \}$. Next we derive the regularity
of $\bar u_{T_0}$ by applying the weighted $W^{2,p}$ estimates in Theorem 3.1. Its regularity
will depend on the relationship between $\ds\f{p_0}{\mu}$ and $Q_0$. Precisely,

\vskip 0.1 true cm

{\bf Case A.} $1<\ds\f{p_0}{\mu}<Q_0$

We have from (3.3) and (3.4) that for any $1< p \le \ds\f{p_0}{\mu}$
$$\bar u_{T_0}\in C\big([-T_0, 0], L^p(\R^n)\cap L^{p_2}(\R^n)\big), \quad \p_t \bar u_{T_0} \in C\big([-T_0, 0], L^p(\R^n)\big)  \cap L^{p_2}(G_{T_0}), \eqno{( 4.5 )} $$ where  $ \ds\f{1}{p_2}=\f{\mu}{p_0}-\f{1}{Q_0}$.

\vskip 0.2 true cm

{\bf Case B.}  $\ds\f{p_0}{\mu}>Q_0$

We have from (3.3) and (3.5) that for any $1< p \le \ds\f{p_0}{\mu}$
$$\bar u_{T_0}\in C\big([-T_0, 0], L^p(\R^n)\cap L^{\infty}(\R^n) \big), \quad \p_t \bar u_{T_0} \in C\big([-T_0, 0], L^p(\R^n)\big)  \cap L^{\infty}(G_{T_0}), $$
which means that for any $1< q < \infty$
$$\bar u_{T_0}\in C\big([-T_0, 0], L^q(\R^n)\cap L^{\infty}(\R^n)\big), \quad \p_t \bar u_{T_0} \in C\big([-T_0, 0], L^q(\R^n)\big)\cap L^{\infty}(G_{T_0}). \eqno{ (4.6)} $$

\vskip 0.2 true cm

{\bf Case C.}  $\ds\f{p_0}{\mu}=Q_0$

We have from (3.3) and (3.6) that for any $1< p < \infty$
$$\bar u_{T_0},  \p_t \bar u_{T_0} \in C\big([-T_0, 0], L^p(\R^n)\big),$$
especially, $$\bar u_{T_0}\in C\big([-T_0, 0], L^{p_0}(\R^n)\big), \quad \p_t \bar u_{T_0} \in C\big([-T_0, 0], L^{p_1}(\R^n)\big). \eqno{( 4.7)}$$

Notice  that for $1< \mu<p_0 $, if $Q_0 \le \ds\f{p_0}{\mu-1}$, there hold $p_2 \ge p_0$ and $\ds\f{p_0}{\mu} > p_1$;
however, if $Q_0 > \ds\f{p_0}{\mu-1}$, then $\ds\f{p_0}{\mu} < p_2<p_0$ and thus  we cannot get the $L^{p_0}$ (only $L^{p_2}$) regularity of $\bar u_{T_0}$ with respect to the variable $x$ from (4.5), which will lead to the difficulty in closing the $L^{p_0}-$ regularity
argument on the solution $u$ to the nonlinear problem (4.1).
Consequently, we will restrict our consideration on the cases for $0 \le \mu \le 1$, or for $1 < \mu <p_0$ with $Q_0 \le \ds\f{p_0}{\mu-1}$.

Collecting (4.5)-(4.7) yields that when $0 \le \mu \le 1$, or when $1 < \mu <p_0$ with $Q_0 \le \ds\f{p_0}{\mu-1}$,
$$\bar u_{T_0}\in C\big([-T_0, 0], L^{p_0}(\R^n)\big), \quad \p_t \bar u_{T_0} \in C\big([-T_0, 0],
L^{p_1}(\R^n)\big). \eqno{(4.8)} $$

In order to solve (4.1) in $\{t\le 0\}$, we only require to consider the following problem
\begin{equation*}
\left\{
\begin{aligned} &\p_t^2 v_{T_0}-t^{2l-1}\Delta v_{T_0}=\chi_{T_0}(t)\bigl(f(t,x,\bar u+\bar u_{T_0}+v_{T_0})
-f(t,x,\bar u)\bigr),\ \  (t,x)\in (-\infty, 0] \times\Bbb R^n,\\
&v_{T_0}(0,x)=0.
\end{aligned}
\right.\tag{4.9}
\end{equation*}
Indeed, if (4.9) is solved, then $u_{T_0}=\bar u+\bar u_{T_0}+v_{T_0}$ is a solution to the following problem
\begin{equation*}
\left\{
\begin{aligned} &\p_t^2 u_{T_0}-t^{2l-1}\Delta u_{T_0}=\chi_{T_0}(t)f(t,x,u_{T_0}),\quad (t,x)\in (-\infty, 0] \times\Bbb R^n,\\
&u_{T_0}(0,x)=u_0(x),
\end{aligned}
\right.
\end{equation*}
namely, we solve the problem (4.1) for $(t,x)\in [-T_0, 0]\times\Bbb R^n$.

We will use the following iteration scheme to study the problem (4.9)
\begin{equation*}
\left\{
\begin{aligned} &\p_t^2 v^{k+1}_{T_0}-t^{2l-1}\Delta v^{k+1}_{T_0}=\chi_{T_0}(t)\bigl(f(t,x,\bar u+\bar u_{T_0}+v^k_{T_0})
-f(t,x,\bar u)\bigr),\quad (t,x)\in  (-\infty, 0] \times\Bbb R^n,\\
&v^{k+1}_{T_0}(0,x)=0,
\end{aligned}
\right.\tag{4.10}
\end{equation*}
where $v^0_{T_0}(t,x)\equiv 0$.

For $k=0$,  $v_{T_0}^1$ is the solution to the following problem
\begin{equation*}
\left\{
\begin{aligned} &\p_t^2 v^{1}_{T_0}-t^{2l-1}\Delta v^{1}_{T_0}=\chi_{T_0}(t)\bigl(f(t,x,\bar u+\bar u_{T_0})
-f(t,x,\bar u)\bigr),\quad (t,x)\in  (-\infty, 0] \times\Bbb R^n,\\
&v^{1}_{T_0}(0,x)=0.
\end{aligned}
\right.\tag{4.11}
\end{equation*}
Denote
$$ I_0\equiv\f 1 t \chi_{T_0}(t) \big( f(t,x, \bar{u}+\bar{u}_{T_0})-f(t, x, \bar{u})\big).$$
Since  $\bar{u}_{T_0}(t,x)=  \int_0^t \p_\tau\bar{u}_{T_0}(\tau, x) d\tau $,  we have from
the condition (1.2) that
\begin{equation*}
|I_0| \le \left\{
\begin{aligned}
 &  C \chi_{T_0}(t)\Big| \ds\f 1 t\int_0^t \p_\tau\bar{u}_{T_0}(\tau, x) d\tau\Big|, \qquad
\qquad \qquad \qquad \qquad \quad \text{ when } 0\le \mu \le 1,\\
& C \chi_{T_0}(t)(1+ |\bar{u}|^{\mu-1} +|\bar{u}_{T_0}|^{\mu-1} ) \Big| \ds\f 1 t\int_0^t \p_\tau\bar{u}_{T_0}(\tau, x) d\tau\Big|,\quad
\text{ when } 1< \mu <p_0.
\end{aligned}
\right.
\end{equation*}
Now we analyze the regularity of $I_0$.
First consider the case $1< \mu <p_0$ with $1< \ds\f {p_0} { \mu} <Q_0 \le \ds\f {p_0} { \mu-1}$. In this case,    $\bar{u} \in L^{p_0}(G_{T_0}),$ $ \bar{u}_{T_0}, \p_t \bar{u}_{T_0} \in L^{p_2}(G_{T_0})$, which together with H\"older's inequality and Hardy's inequality derives that
$$ I_0 \in L^{p_2}(\R^{n+1}_-) +L^{ {p_2}/{\mu}}(\R^{n+1}_-)+ L^{r_1}(\R^{n+1}_-),\eqno{(4.12)} $$
where and below $\ds\f{1}{r_1}=\f{\mu-1}{p_0}+\f 1 {p_2}=\f{2\mu-1}{p_0}-\f{1}{Q_0}$.

Since $I_0$ has compact support with respect to variables $t$ and $x$, and in this case, $\min\{p_2, \ds\f {p_2}{\mu}, r_1\}$ $=r_1,$
we have that for any $1 < \gamma\le r_1$,
$$
I_0\in L^{\gamma}(\R^{n+1}_-),   \quad \text{ when } 1< \mu <p_0  \ \ \text{ and } 1< \ds\f {p_0} { \mu} <Q_0 \le \ds\f {p_0} { \mu-1}.\eqno{(4.13)}
$$
Similarly, we also have that
\begin{equation*}
I_0 \in \left\{
\begin{aligned}
 & L^{p}(\R^{n+1}_-),  \qquad \qquad \qquad \text{ for any }  1< p \le \ds p_2, \quad \text{ if $0< \mu \le 1$ and
$1<\f{p_0}{\mu}<Q_0$}, \\
&L^{q}(\R^{n+1}_-) \cap L^{\infty}(\R^{n+1}_-), \ \ \text{  for any }  1< q <\infty, \quad \text{ if $0< \mu \le 1$ and
$\ds\f{p_0}{\mu}>Q_0$}, \\
&L^\delta(\R^{n+1}_-),  \qquad \qquad \qquad \text{ for any }  1< \delta \le \ds\f{p_0}{\mu-1}, \quad \text{ if $1< \mu <p_0$
and $\ds\f {p_0}{\mu} >Q_0$}, \\
&L^\eta(\R^{n+1}_-), \qquad \qquad \qquad \text{ for any }  1< \eta \le \ds\f{p_0}{\mu}, \quad \quad \text{ if $\ds\f {p_0}{\mu}= Q_0$}.
\end{aligned}
\right.\tag{4.14}
\end{equation*}
Based on (4.13)-(4.14), we then have from (4.11) and Theorem 3.1 with $m=2l-1$ and $\nu=1$ that

\vskip 0.2 true cm

(i) For $0< \mu \le 1,$
 \begin{align*}
 v_{T_0}^1\in \hskip 16cm   \\
\left\{
\begin{aligned}
 & C([-T_0, 0], L^q(\R^n)\cap L^\infty(\R^n))\quad\text{for $1<q< \infty$},  \quad  \text{ if $\ds\f {p_0}{\mu} > Q_0$,}
\text{ or if $\ds\f {p_0}{\mu} < Q_0$ and $p_2 >Q_1$}, \\
& C\big([-T_0, 0], L^p(\R^n)   \cap L^{\theta}(\R^n)\big)\quad\text{for $1< p \le p_2$ and $\ds\f 1 {\theta}\equiv\ds \f 1 {p_2} -\f 1{Q_1}$,}\quad \text{if $\ds\f {p_0}{\mu} < Q_0$ and $ p_2 <Q_1$},\\
\end{aligned}
\right.
 \end{align*}

and
\begin{align*}
 \p_t v_{T_0}^1 \in \hskip 15cm   \\
\left\{
\begin{aligned}
 & C([-T_0, 0], L^q(\R^n)) \cap L^\infty(G_{T_0})\quad\text{for $1<q< \infty$},  \text{ if $\ds\f {p_0}{\mu} > Q_0$,
or if $\ds\f {p_0}{\mu} < Q_0$ and $p_2 >Q_1$}, \\
& C\big([-T_0, 0], L^p(\R^n)\big)   \cap L^{\theta}(G_{T_0})\quad\text{for $1< p \le p_2$},    \quad  \text{ if $\ds\f {p_0}{\mu} < Q_0
$ and $1< p_2 <Q_1$.}
\end{aligned}
\right.
\end{align*}

(ii) For $1<\mu <p_0,$
\begin{align*}
v_{T_0}^1\in
\left\{
\begin{aligned}
 & C([-T_0, 0], L^p(\R^n) \cap L^\infty(\R^n))\quad\text{ for $1<p< \infty$}, \  \text{ if $\ds\f {p_0}{\mu} > Q_0$,}\\
&\qquad \qquad \qquad \qquad \qquad \qquad \qquad \qquad \text{ or if $ \ds\f {p_0}{\mu}<Q_0 \le \ds\f {p_0} { \mu-1}$ and $r_1>Q_1$}, \\
& C([-T_0, 0], L^\gamma(\R^n)   \cap L^{\Theta}(\R^n)) \text{for $1< \gamma \le r_1$,  $\ds\f 1 \Theta\equiv\f 1 r_1 - \f 1 {Q_1}$},\\
&\qquad \qquad \qquad \qquad \qquad \qquad \qquad \quad  \text{if $\ds \f {p_0}{\mu}<Q_0\le \ds\f {p_0} { \mu-1}$ and $ r_1 <Q_1$,}\\
\end{aligned}
\right.
\end{align*}

and
\begin{align*}
\p_t v_{T_0}^1 \in
 \left\{
\begin{aligned}
 & C([-T_0, 0], L^p(\R^n))  \cap L^\infty(G_{T_0})\  \text{for $1<p< \infty$,}  \ \text{ if $\ds\f {p_0}{\mu} > Q_0$,}\\
&\qquad\qquad \qquad \qquad \qquad \qquad \qquad  \text{or if $1< \ds\f {p_0}{\mu}<Q_0\le \ds\f {p_0} { \mu-1}$, $r_1>Q_1$}, \\
& C([-T_0, 0], L^\gamma(\R^n))   \cap L^{\Theta}(G_{T_0})  \ \text{for $1< \gamma \le r_1$},\\
&\qquad \qquad \qquad \qquad \qquad \qquad \text{ if $1< \ds\f {p_0}{\mu}<Q_0\le \ds\f {p_0} { \mu-1}$
and $1< r_1 <Q_1$}.
\end{aligned}
\right.
\end{align*}

Moreover, when  $\ds\f {p_0}{\mu} =Q_0$, or when $0< \mu \le 1$ with $\ds\f {p_0}{\mu} < Q_0$ and $ p_2 =Q_1$,
or when $1<\mu <p_0$ with $\ds \f {p_0}{\mu}<Q_0\le \ds\f {p_0} { \mu-1}$ and $ r_1 =Q_1$, we have that
$$ v_{T_0}^1, \p_t v_{T_0}^1 \in C([-T_0, 0], L^q(\R^n)) \quad \text{ for any $1<q <\infty$}. $$
Therefore, collecting all the above analysis in (i)-(ii)  yields
$$v_{T_0}^1\in C\Big((-\infty, 0], L^{p_0}(\Bbb R^n)   \Big), \ \ \p_t  v_{T_0}^1\in C\Big((-\infty, 0], L^{p_1}(\Bbb R^n)\Big), \eqno{(4.15)}$$
when $0 < \mu \le 1$, or when $1 < \mu <p_0$ and $Q_0 \le \ds\f{p_0}{\mu-1}$.

To study the problem (4.10), we denote  by
$$I_k\equiv\f 1 t \chi_{T_0}(t) \Big( f(t,x, \bar{u}+\bar{u}_{T_0}+ v_{T_0}^k)-f(t, x, \bar{u}) \Big).$$

We now derive the regularities of $I_k$ under  different restrictions on $\mu$, $l$ and $p_0$
in (1.2)-(1.3), whose classifications are based on  the distinguished regularities of $v_{T_0}$ and $\p_t v_{T_0}$.

Define the set $\mathcal M_1$
\begin{align*}
\mathcal M_1=\{g(t,x): & \ g(t,x)\in C([-T_0, 0], L^q(\R^n)\cap L^\infty(\R^n)), \\
&\p_t g(t,x) \in C([-T_0, 0], L^q(\R^n) ) \cap L^\infty(G_{T_0}), \quad
\||g|\|_1+ \t\||\p_tg|\t\|_1\le 2\},\\
& \text{ for any } 1<q<\infty, \end{align*}
where $$\||g|\|_1\equiv\|g\|_{C([-T_0, 0], L^q(\R^n))}+ \|g\|_{C([-T_0, 0],  L^\infty(\R^n))}$$ and
$$\t\||\p_t g|\t\|_1\equiv\|\p_tg\|_{C([-T_0, 0], L^q(\R^n))}+\|\p_tg\|_{L^{\infty}(G_{T_0})}.$$ In particular,
for $g \in \mathcal M_1$, one has $g\in C([-T_0, 0], L^{p_0}(\R^n)) \cap C^1([-T_0, 0], L^{p_1}(\R^n)).$

\vskip .1cm

{\bf Case 1.} $1<\mu <p_0$ and $Q_0<\ds\f{p_0}{\mu}$

In this case, for $v_{T_0}^k \in \mathcal M_1,$ we have from (1.2) that
$$|I_k |  \le  C \chi_{T_0}(t)(1+ |\bar{u}|^{\mu-1} +|\bar{u}_{T_0}|^{\mu-1} +|v_{T_0}^k|^{\mu-1}   ) \Big( | \f 1 t\int_0^t \p_\tau\bar{u}_{T_0}(\tau, x) d\tau| + | \f 1 t\int_0^t \p_\tau v_{T_0}^k(\tau, x) d\tau|\Big),$$
which, together with  Hardy's inequality and the regularities of $\bar{u}, \bar{u}_{T_0}, v_{T_0}^k, \p_t \bar{u}_{T_0}$ and $\p_t v_{T_0}^k$ in $G_{T_0}$ (that is, $\bar{u} \in L^{p_0}(G_{T_0})$ and $\bar{u}_{T_0},$ $\p_t \bar{u}_{T_0}$, $v_{T_0}^k$, $\p_t v_{T_0}^k \in L^{\infty}(G_{T_0})$), yields
$$I_k \in L^\infty(\R^{n+1}_-) + L^{\f {p_0}{\mu-1}}(\R^{n+1}_-).$$
Since $I_k$ has  the compact support with respect to the variables $t$ and $x$, we have
$$I_k \in L^{\delta}(\R^{n+1}_-) \text{ for any $1 < \delta\le \ds\f {p_0}{\mu-1}$}.$$
Noting that in this case, $\ds\f {p_0}{\mu -1} >\ds\f {p_0}{\mu } >Q_0>Q_1$ hold, then we have from (4.10),  (3.5) and (3.36)-(3.37) that $v_{T_0}^{k+1} \in \mathcal M_1$
for small $T_0>0$.

\vskip .2cm

{\bf Case 2.} $1<\mu <p_0,$ $\ds\f {p_0}{\mu} < Q_0 \le \ds\f {p_0}{\mu -1} $ and $r_1 >Q_1$

In this case, for any $v_{T_0}^k \in \mathcal M_1,$
$$\bar{u} \in L^{p_0}(G_{T_0}), \ \bar{u}_{T_0}\in L^{p_2}(G_{T_0}), \ \p_t \bar{u}_{T_0}\in L^{p_2}(G_{T_0}), \ v_{T_0}^k\in L^{\infty}(G_{T_0}),  \ \p_t v_{T_0}^k \in L^{\infty}(G_{T_0}),$$
and then as in Case 1,  we have in $G_{T_0}$,
$$ I_k \in L^{p_2} + L^\infty+ L^{r_1}+ L^{p_2/\mu } + L^{p_2/ (\mu -1) } +  L^{p_0/ (\mu -1) }.$$
Noting
$$\min\{p_2,  r_1,\f {p_2} {\mu}, \f { p_2} { \mu -1}, \f { p_0} {\mu -1} \}= r_1$$
and $I_k$  has  the compact support on the variables $t$ and  $x$, then
$$I_k \in L^{\gamma}(\R^{n+1}_-) \qquad \text{ for any $1 < \gamma\le r_1$.}$$
Due to our assumption in this case, $r_1 > Q_1$ holds. Thus it follows from (3.36)-(3.37) and (3.5) in Theorem 3.1, and (4.10)  that $v_{T_0}^{k+1} \in \mathcal M_1$
for small $T_0>0$.

\vskip .2cm

{\bf Case 3.}  $0<\mu\le 1,$ $Q_0>\ds\f {p_0}{\mu}$ and $p_2 >Q_1$

In this case, for $v_{T_0}^k \in \mathcal M_1,$  the condition (1.2) implies that
$$  |I_k| \le  C \chi_{T_0}(t)\Big( \big| \f 1 t\int_0^t \p_\tau\bar{u}_{T_0}(\tau, x) d\tau\big|
 + \big| \f 1 t\int_0^t \p_\tau v_{T_0}^k(\tau, x) d\tau\big| \Big),$$
which yields
$$ I_k \in L^\infty(\R^{n+1}_-) + L^{p_2}(\R^{n+1}_-) \subset L^{p_2}(\R^{n+1}_-) $$
and then
$$I_k \in L^{p}(\R^{n+1}_-)\quad \text{ for any $1 < p \le p_2$.}$$
By $p_2 >Q_1$ from our assumption, then estimates   (3.5), (3.36)-(3.37) in Theorem 3.1 and equation  (4.10) that $v_{T_0}^{m+1} \in \mathcal M_1$
for small $T_0>0$.

\vskip .2cm

{\bf Case 4.} $0<\mu\le 1$ and $Q_0<\ds\f {p_0}{\mu}$

For $v_{T_0}^k \in \mathcal M_1,$  we have
$$I_k \in L^{p}(\R^{n+1}_-) \cap L^{\infty}(\R^{n+1}_-)\quad \text{ for any $1 < p<\infty$},$$
thus, as in Case 3, $v_{T_0}^{k+1} \in \mathcal M_1$ for small $T_0>0$.

\vskip .2cm

In order to get  suitable regularities of $I_k$ for some other left cases of $\mu, l$ and $p_0$,
we now define the second set $\mathcal M_2$ as follows
\begin{align*}  \mathcal M_2=\{g(t,x): \ & g\in C([-T_0, 0], L^p(\R^n) \cap L^{\theta}(\R^n)),
\p_t g\in C([-T_0, 0], L^p(\R^n) ) \cap L^{\theta}(G_{T_0}), \\
& \hskip 5cm \||g|\|_2+ \t\||\p_tg|\t\|_2\le 2\},\end{align*}
where   $$\||g|\|_2\equiv\|g\|_{C([-T_0, 0], L^p(\R^n))}
+\|g\|_{C([-T_0, 0], L^\th(\R^n))}$$ and $$\t\||\p_tg|\t\|_2\equiv\|\p_tg\|_{C([-T_0, 0], L^p(\R^n))}
+\|\p_tg\|_{L^{\theta}(G_{T_0})},$$  here $1<p\le p_2$, $\ds\f{1}{\th}=\f{1}{p_2}-\f{1}{Q_1}$.
We now analyze the regularity of $I_k$ when $v_{T_0}^k \in \mathcal M_2$.

\vskip .2cm

{\bf Case 5.} $0<\mu \le1,$ $Q_0>\ds\f{p_0}{\mu}$ and $Q_1>p_2$

For $v_{T_0}^k \in \mathcal M_2,$   we have from Hardy's inequality and the regularity of  $\p_t\bar{u}_{T_0} \in L^{p_2}(G_{T_0})$ and $\p_t v_{T_0}^k \in L^{\theta}(G_{T_0})$  that
$$I_k \in  L^{p_2}(\R^{n+1}_-)+ L^{\theta}(\R^{n+1}_-)\subset L^{p_2}(\R^{n+1}_-)$$
and further
$$I_k \in  L^{p}(\R^{n+1}_-) \quad\text{ for any $1 < p \le p_2$.}$$
Due to $1< p_2 < Q_1$ from our assumption in this case, we have from the estimates (3.36)-(3.37) and (3.4) in
Theorem 3.1 together with (4.10)  that $v_{T_0}^{k+1} \in \mathcal M_2$ for small $T_0>0$. Obviously,  under the assumption in Case 5, if $g \in \mathcal M_2$, then
$$g\in C([-T_0, 0], L^{p_0}(\R^n)) \cap C^1([-T_0, 0], L^{p_1}(\R^n)).$$

As before, we require to
define the third set $\mathcal M_3$
\begin{align*}  \mathcal M_3=\{g(t,x): \ & g\in C([-T_0, 0], L^\gamma(\R^n) \cap L^{\Theta}(\R^n)),  \p_t g
\in C([-T_0, 0], L^\gamma(\R^n) ) \cap L^{\Theta}(G_{T_0}),  \\
&\hskip 5cm \||g|\|_3+ \t\||\p_tg|\t\|_3\le 2\},  \end{align*}
where  $$\||g|\|_3\equiv\|g\|_{C([-T_0, 0], L^\gamma(\R^n))}
+\|g\|_{C([-T_0, 0], L^\Theta(\R^n))}$$ and $$\t\||\p_tg|\t\|_3\equiv\|\p_tg\|_{C([-T_0, 0], L^\gamma(\R^n))}
+\|\p_tg\|_{L^{\Theta}(G_{T_0})},$$ here $1<\gamma\le r_1, \ds\f{1}{\Theta}\equiv\f{1}{r_1}-\f{1}{Q_1}.$ We now derive the regularity of $I_k$ in the following case.

\vskip .2cm

{\bf Case 6.} $1<\mu <p_0,$ $\ds\f{p_0}{\mu} < Q_0 \le \f{p_0}{\mu -1}$ and $Q_1>r_1$

At this moment, we have $\Theta > p_0$ and $ r_1 > p_1$. Hence,  in Case 6, if $g \in \mathcal M_3$, then
$$g\in C([-T_0, 0], L^{p_0}(\R^n)) \cap C^1([-T_0, 0], L^{p_1}(\R^n)),$$
and for $v_{T_0}^k \in \mathcal M_3$,
$$I_k \in L^{p_2} + L^{\Theta} + L^{r_1} + L^{ {p_2}/{\mu}}+L^{\eta_1}+L^{\eta_2} +L^{\eta_3} +L^{\Theta/\mu},$$
where
$$
\ds\f 1 {\eta_1}\equiv\f {\mu-1}{\Theta}+ \f 1{p_2}, \quad \f 1 {\eta_2}\equiv\f {\mu-1}{p_0}+ \f 1{\Theta},
\quad \f 1 {\eta_3}\equiv\f {\mu-1}{p_2}+ \f 1{\Theta}.
$$
Noting that
$$\min\{p_2, \Theta, \f {p_2} { \mu}, \f {\Theta} { \mu}, r_1, \eta_1, \eta_2, \eta_3\}=r_1,$$
then
$$I_k \in L^{\eta}(\R^{n+1}_-)\quad \text{ for any $1 < \eta\le r_1$.}$$
This  together with (3.4) and (3.36)-(3.37) in Theorem 3.1  yields $v_{T_0}^{m+1} \in \mathcal M_3$ for small $T_0>0$.

Finally, we treat the term $I_k$ for the left cases of $\mu, l$ and $p_0$. To this end, as before,
we define the fourth set $\mathcal M_4$
$$
\mathcal M_4=\{g(t,x): \  g,  \p_t g\in C([-T_0, 0], L^{p_0}(\R^n) ),  \||g|\|_4+ \t\||\p_tg|\t\|_4\le 2\},$$
where $\||g|\|_4=\|g\|_{C([-T_0, 0], L^{p_0}(\R^n) )}$ and $\t\||\p_tg|\t\|_4=\|\p_tg\|_{C([-T_0, 0], L^{p_0}(\R^n))}$.
We now distinguish the following three cases:

\vskip 0.2 true cm
{\bf Case 7.}  $0< \mu < p_0$ and $Q_0=\ds\f{p_0}{\mu}$

\vskip 0.1 true cm

{\bf Case 8.}  $1< \mu < p_0$, $\ds\f{p_0}{\mu}<Q_0\le \ds\f{p_0}{\mu}$ and $Q_1=r_1$

\vskip 0.1 true cm
{\bf Case 9.}  $0< \mu\le 1$, $Q_0>\ds\f{p_0}{\mu}$ and $Q_1= p_2$

\vskip 0.1 true cm

By the analogous analysis as in the previous cases, we can obtain $v_{T_0}^{k+1} \in \mathcal M_4$
for the Case 7 -Case 9 when $T_0>0$ is small.

\vskip .2cm

Next we show the convergence of the sequence $\{v_{T_0}^k\}$ in the corresponding spaces $\mathcal M_i$
$(1\le i\le 4)$. Set $\dot v_{T_0}^{k+1}=v_{T_0}^{k+1}-v_{T_0}^k$, then it follows from (4.10) that
\begin{equation*}
\left\{
\begin{aligned} &\p_t^2\dot v_{T_0}^{k+1}-t^{2l-1}\Delta \dot v_{T_0}^{k+1}=\chi_{T_0}(t)
\bigl(f(t,x,\bar u+\bar u_{T_0}+v^k_{T_0})
-f(t,x,\bar u+\bar u_{T_0}+v^{k-1}_{T_0})\bigr),\\
&\dot v_{T}^{k+1}(0,x)=0.
\end{aligned}
\right.
\end{equation*}
As in Case 1-Case 9, in the related spaces $\mathcal M_j$ $(1\le j\le 4)$,
we can obtain by (3.35)-(3.37) separately
$$\||\dot v_{T_0}^{k+1}|\|_j+\||\p_t\dot v_{T_0}^{k+1}|\|_j  \le
C(T_0)\||\dot v_{T_0}^{k}|\|_j,\eqno{(4.16)} $$
where $C(T_0)\le\ds\f12$ for small $T_0>0$. Therefore, one can derive from (4.16) that
for small fixed  $T_0>0$, there exists a function $v_T(t,x)$ such that
$v_{T_0}^{k+1}\to v_{T_0}(t,x)$ in $C( [-T_0, 0], L^{p_0}(\R^n))  \cap C^1( [-T_0, 0], L^{p_1}(\R^n))$
and $v_{T_0}\in C( [-T_0, 0], L^{p_0}(\R^n))\cap C^1( [-T_0, 0], L^{p_1}(\R^n)) $  solves
(4.9). Therefore, $u=\bar u+{\bar u}_{T_0}+v_{T_0}$ solves the problem (4.1) for $(t,x)\in [-T_0, 0]\times\Bbb R^n$,
and we complete the proof of Theorem 4.1  \hfill $\square$

\vskip .3cm

{\bf Remark 4.1.} {\it If the assumptions of $0\le\mu\le 1$ or $1<\mu<p_0$ with $Q_0\le\ds\f{p_0}{\mu-1}$
are violated, then we cannot use the standard iteration scheme (4.10) to derive the existence of the solution
to the problem (4.9). Indeed, for example,
when $1< \mu <p_0$ with $ Q_0 > \ds\f {p_0}{\mu -1}$ and $Q_1 > \ds\f { p_1}{\mu-1}$, if  the
iteration scheme (4.10) does work, then it follows from Theorem 3.1 that for $v_{T_0}^{k} \in L^{q_{k}}(\R^{n+1}_-)$
$$ v_{T_0}^{k+1} \in L^{q_{k+1}}(\R^{n+1}_-), \qquad \text{where $\ds\f {1}{q_{k+1}}= \f {\mu}{q_k} -\f 1 {Q_1}$.} $$
Noticing that one has for $k \ge 1$
$$\ds\f 1 {q_{k+1}} - \f 1 {q_{k}} =\mu^{k-1} \big( \f 1 {q_2} - \f 1{q_1} \big)=\mu^k \big( \f{\mu-1}{p_1} -\f {1}{Q_1} \big),
$$
thus, we have
\begin{align*}
\f 1 {q_k}=&  \f{\mu^k}{p_1} -\f {\mu^k - 1}  { ( \mu -1)  Q_1 } =\mu^k \Big(  \f 1 {p_1} - \f 1 {(\mu-1)Q}\Big)+ \f 1 {(\mu-1)Q_1} \\
\ge & \mu^k \Big(  \f 1 {p_1} - \f 1 {(\mu-1)Q}\Big) \rightarrow +\infty \quad \text{ as } m\rightarrow +\infty.
\end{align*}
Therefore, it seems that it is difficult for us to derive the uniform lower bound of $q_k>1$ for any $k\in\Bbb N$.
This means that the standard iteration scheme (4.10) only works for finite steps.}

\vskip 0.4 true cm

\section{ Solvability of (1.1) in the degenerate hyperbolic region $\{t\ge 0\}$ }
Based on Theorem 4.1 in $\S 4$, we consider the  problem (1.1) in the degenerate hyperbolic region $\{t\ge 0\}$
\begin{equation*}
\left\{
\begin{aligned} &\p_t^2 w-t^{2l-1}\Delta w=f(t,x,w),\qquad (t,x)\in [0, \infty)\times\Bbb R^n,\\
&w(0,x)=\varphi(x),\\
&  \p_t w(0,x)=\psi(x),
\end{aligned}
\right.\tag{5.1}
\end{equation*}
where $\varphi(x)\in H^s(\R^n)\subset L^{p_0}(\R^n)$, $0 \le s <\ds\f n 2$, and $\psi(x)\equiv\p_t u(0,x)\in L^{p_1}(\Bbb R^n)$,
here $\p_tu(0,x)$  comes from
Theorem 4.1 in the degenerate elliptic part $\{t \le 0\}$ of the problem (1.1). Our main result in this section is
\vskip .2cm

{\bf Theorem 5.1.} {\it Under the conditions (1.2)-(1.3), there exists a small constant $T_0>0$ such that the problem (5.1)
has a unique local solution $w(t,x)\in C([0, T], L^{p_0}(\Bbb R^n))$ when $0 \le \mu \le 1$, or when $1 < \mu <p_0$
with $Q_0 \le \ds\f{p_0}{\mu-1}$.
}

\vskip .1cm

{\bf Proof.}  We first consider the following linear problem
\begin{equation*}
\left\{
\begin{aligned} &\p_t^2 w_1-t^{2l-1}\Delta w_1=0,\quad (t,x)\in [0, \infty)\times\Bbb R^n,\\
&w_1(0,x)=\varphi(x), \\
&  \p_t w_1(0,x)=\psi(x).
\end{aligned}
\right.\tag{5.2}
\end{equation*}
By [26] or Lemma 2.2 in [19], we know that the problem (5.2) has a unique  solution $w_1$, which    can be expressed as follows
$$w_1(t,x) =V_1(t, D_x) \varphi(x) +V_2(t, D_x) \psi(x),  \eqno{( 5.3)}$$
 where the pseudo-differential  operator  $V_j(t, D_x)$ has the  symbol $V_j(t, |\xi|)$ for $j=1,2,$
\begin{equation*}
\left\{
\begin{aligned} &V_1(t,
|\xi|)=e^{-\f{z}{2}}\Phi(\ds\f{2l-1}{2(2l+1)},\f{2l-1}{2l+1};z),\\
&V_2(t,|\xi|)=te^{-\f{z}{2}}\Phi(\ds\f{2l+3}{2(2l+1)},\f{2l+3}{2l+1};z),
\end{aligned}
\right.\tag{5.4}
\end{equation*}
here the confluent hypergeometric functions $\Phi(\ds\f{2l-1}{2(2l+1)},\f{2l-1}{2l+1};z)$ and
$\Phi(\ds\f{2l+3}{2(2l+1)},\f{2l+3}{2l+1};z)$ are analytic functions of
$z$ with  $z=\ds\f{4i}{2l+1}t^{\f{2l+1}{2}}|\xi|$. Moreover, for
sufficiently large $|z|$,
$$|\Phi(\f{2l-1}{2(2l+1)},\f{2l-1}{2l+1};z)|\le C|z|^{-\f{2l-1}{2(2l+1)}}
\big(1+O\big(|z|^{-1}\big)\big), $$
and
$$ |\Phi(\f{2l+3}{2(2l+1)},\f{2l+3}{2l+1};z)|\le C|z|^{-\f{2l+3}{2(2l+1)}}
\big(1+O\big(|z|^{-1}\big)\big). \eqno{( 5.5)} $$

In addition, by  Lemma 3.2 in [19] and Sobolev's embedding theorem, one has
$$\|V_1(t, D_x) \varphi(x)\|_{L^{p_0}(\R^n)} \le C \|V_1(t, D_x) \varphi(x)\|_{H^s(\R^n)}\le C \|\varphi\|_{H^s(\R^n)}\eqno{( 5.6)} $$
and
$$\|\p_t V_1(t, D_x) \varphi(x)\|_{L^{q_0}(\R^n)} \le C  \|\p_t V_1(t, D_x) \varphi(x)\|_{H^{s-\f {2l+3}{2(2l+1)}}(\R^n)} \le C \|\varphi\|_{H^s(\R^n)}, \eqno{( 5.7)}$$
where $\ds\f 1 {q_0} = \f 1 2 - \f{1}{n}\biggl(s- \f {2l+3}{2(2l+1)}\biggr).$

On the other hand, it follows from the analytic property of $\Phi(\ds\f{2l+3}{2(2l+1)},\f{2l+3}{2l+1};z)$
on the variable $z$ and the estimate (5.5) that there exists a constant $C>0$ such that
$$|z|^{\f 2 {2l+1}} |V_2(t, |\xi|)| \le C t$$
and further
$$|V_2(t, |\xi|)| \le C |\xi|^{-\f 2 {2l+1} }.  \eqno{( 5.8)} $$
It is noted that $\psi(x)=\p_t u(0,x) \in L^{p_1}(\R^n)$ and
$$1< p_1 <p_0 <+\infty, \ \ \f 1 {p_1} -\f 1 {p_0}= \f 2 {n(2l+1)}.$$
This, together with (5.8) and Hardy-Littlewood-Sobolev's inequality,   yields
$$ \|V_2(t, D_x) \psi(x)\|_{L^{p_0}(\R^n)} \le C \| \psi\|_{L^{p_1}(\R^n)}. \eqno{( 5.9)} $$
Therefore, by (5.3), (5.6) and (5.9), we obtain $$w_1(t,x) \in C\big([0,T], L^{p_0}(\R^n)\big).\eqno{( 5.10)} $$

Set $v(t,x)= w(t,x)- w_1(t,x)$, then it follows from (5.1) that
\begin{equation*}
\left\{
\begin{aligned} &\p_t^2 v-t^{2l-1}\Delta v=f(t,x,v+w_1),\qquad (t,x)\in [0, \infty)\times\Bbb R^n,\\
&v(0,x)=0,\qquad \p_t v(0,x)=0.
\end{aligned}
\right.\tag{5.11}
\end{equation*}
For suitably chosen constant $T>0,$ we define a set $\mathcal M$:
$$ \mathcal M=\{ v \in C\big([0, T], L^{p_0}(\R^n)\big): \sup_{t \in [0,T]} \|v(t, \cdot)\|_{L^{p_0}(\R^n)} \le 2\}. $$
And then we define a mapping $\mathcal T$ as follows
$$\mathcal T(v)(t, x)= \int_0^t \Big(V_2(t, D_x) V_1(\tau, D_x)- V_1(t, D_x) V_2(\tau, D_x)\Big) f(\tau,x,v+w_1) d\tau.$$
It is easy to verify that $\mathcal T (v)$ satisfies
\begin{equation*}
\left\{
\begin{aligned} &\Big(\p_t^2 -t^{2l-1}\Delta \Big)\mathcal T (v)=f(t,x,v+w_1),\qquad (t,x)\in [0, \infty)\times\Bbb R^n,\\
&\mathcal T (v)(0,x))=0, \\
&  \p_t \mathcal T (v)(0,x)=0.
\end{aligned}
\right.\tag{5.12}
\end{equation*}
By Theorem 3.1 in [26], we have that for any $g \in L^p (\R^n)$ with $1< p< \infty$
$$ \|V_1(t, D_x) g(x)\|_{L^p (\R^n)} \le C \|g\|_{L^p (\R^n)}\eqno{( 5.13)} $$
and
$$ \| V_2(t, D_x) g(x)\|_{L^p (\R^n)} \le C  t \|g\|_{L^p (\R^n)}. \eqno{( 5.14)}$$

Next we show that the mapping $\mathcal T$ maps $\mathcal M$ into itself and the mapping $\mathcal T$ is contractible
for small $T>0$.
If so, then we can get the solvability of (5.11). For this purpose, we will distinguish two cases as follows.

\vskip  .2cm
{\bf Case 1.}  $0\le \mu \le 1:$

\vskip  .2cm

in this case, one has
\begin{align*}
\|\mathcal T(v)(t, \cdot)\|_{L^{p_0}(\R^n)}\le & C  t \int_0^t  \|f(\tau, v+ w_1)\|_{L^{p_0  }(\R^n)} d\tau \quad \text{ (by  (5.13) and (5.14)) } \\
\le & C  t \int_0^t  \|f(\tau, v+ w_1)\|_{L^{p_0 / \mu }(\R^n)} d\tau \quad \text{ (by H\"older's inequality) }\\
\le & C t \int_0^t \bigl(1+ \| v(\tau, \cdot)+ w_1(\tau, \cdot))\|_{L^{ p_0 }(\R^n)}\bigr) d\tau  \quad \text{ (by (1.2)) }\\
\le & C t^2 \Big(1+\sup_{t \in [0, T]} \|v(t, \cdot)\|_{L^{ p_0 }(\R^n)}+ \|w_1(t, \cdot)\|_{L^{ p_0 }(\R^n)} \Big).  \tag {5.15}
\end{align*}

{\bf Case 2.} $1< \mu<p_0$ and $Q_0 \le \ds\f{p_0}{\mu-1}:$

\vskip 0.2 true cm

in this case, one has $p_1 \le\ds\f{p_0}{\mu}<p_0$. Thus, we have
\begin{align*}
\|\mathcal T(v)(t, \cdot)\|_{L^{p_0}(\R^n)}
\le & C  \int_0^t  \|f(\tau, v+ w_1)\|_{L^{ p_1 }(\R^n)} d\tau  \quad \text{ (by  (5.11) and (5.13))}\\
\le & C  \int_0^t  \|f(\tau, v+ w_1)\|_{L^{  p_0 / \mu }(\R^n)} d\tau \quad \text{ (by H\"older's inequality) } \\
\le & C \int_0^t \bigl(1+ \| v(\tau, \cdot)+ w_1(\tau, \cdot))\|_{L^{ p_0 }(\R^n)}\bigr) d\tau \quad \text{ (by (1.2)) } \\
\le & C t \Big(1+\sup_{t \in [0, T]} \|v(t, \cdot)\|_{L^{ p_0 }(\R^n)}+ \|w_1(t, \cdot)\|_{L^{ p_0 }(\R^n)} \Big).\tag {5.16}
\end{align*}
From (5.6), (5.9) and (5.16), we see that for small $T>0$, $\mathcal T(v)(t, x) \in C([0,T], L^{p_0}(\R^n))$ and
$$
\sup_{t \in [0,T]}\| \mathcal T(v)(t, \cdot)\|_{L^{p_0}(\Bbb R^n)} \le   C_{p_0}T
\Big(1+ \sup_{t \in [0,T]}\|  v(t, \cdot)\|_{L^{p_0}(\Bbb R^n)}+\sup_{t \in [0,T]}\|  w_1(t, \cdot)\|_{L^{p_0}(\Bbb R^n)}  \Big)
\le 2,\eqno{(5.17)}
$$
which implies that $\mathcal T$ maps $\mathcal M$ into itself.

In addition, for any $v_1, v_2 \in \mathcal M,$ as argued in (5.16), one can obtain for small $T>0$
$$
\| \mathcal T(v_1)(t, \cdot)-\mathcal T(v_2)(t, \cdot)\|_{L^{p_0}(\Bbb R^n)} \le   \f 1 2
\sup_{t \in [0,T]}\|  v_1(t, \cdot)- v_2(t, \cdot)\|_{L^{p_0}(\Bbb R^n)}.\eqno{(5.18)}
$$
Therefore, by (5.17)-(5.18) and the contraction map principle,
we complete the proof of Theorem 5.1. \hfill $\square$
\vskip 0.2 true cm
{\bf Remark 5.1.} {\it If the initial data $u(0,x)=\vp_0(x)\in L^{\infty}(\Bbb R^n)\cap H^s(\Bbb R^n)$ in (1.1)
for $n=2, 3$ and $s\ge 0$ is given, we can remove all the assumptions in (1.2)-(1.3). Indeed, from the proof of Theorem 4.1,
we can derive that $u(t,x), \p_tu(t,x)\in L^{\infty}([-T_0, 0]\times\Bbb R^n)$ for some fixed
constant $T_0>0$ without the assumptions (1.2)-(1.3), which obviously means $(u(0,x), \p_tu(0,x))=(\vp_0(x), \vp_1(x))\in L^{\infty}(\Bbb R^n)$.
In addition, by [25], one has the following formula
\begin{align*}
w(t,x)=& 2^{2-2\g}\ds\f{\G(2\g)}{\G^2(\g)}\int_0^1v_{\vp_0}(\phi(t)s, x)(1-s^2)^{\g-1}ds \\ & +
2^{2\g}\ds\f{\G(2-2\g)}{\G^2(1-\g)}t\int_0^1v_{\vp_1}(\phi(t)s, x)(1-s^2)^{-\g}ds
\end{align*}
when $w(t,x)$ is a solution
to the equation $\p_t^2w-t^{2l-1}\Delta w=0$ with $(w(0,x), \p_tw(0,x))=(\vp_0(x),$ $\vp_1(x))$, where $\g=\ds\f{2l-1}{2(2l+1)}$,
$\phi(t)=\ds\f{2}{2l+1}t^{\f{2l+1}{2}}$, and $v_{\vp}$ denotes by the solution to
the wave equation $\p_t^2v-\Delta v=0$ with $(v(0,x), \p_tv(0,x))=(\vp(x), 0)$. It is well-known that
$v\in L^{\infty}([0, T]\times\Bbb R^n)$ holds for $n=2,3$ and any $T>0$ when $\vp(x)\in L^{\infty}(\Bbb R^n)$.
Based on this, one easily knows $w(t,x)\in L^{\infty}([0, T]\times\Bbb R^n)$ $(n=2,3)$. Therefore,
by a simpler proof than that in Theorem 5.1, we can locally solve the problem (1.1) with no the assumptions (1.2)-(1.3)
when $n=2,3$, and further $u(t,x)\in L^{\infty}([-T_0, T_0]\times\Bbb R^n)$ can be derived
and (1.1) is solved.
}

\vskip 0.4 true cm

\section{ Proof of Theorem 1.1. }

In this section, based on Theorem 4.1 and Theorem 5.1, we will complete the proof of Theorem 1.1.

{\bf Proof of Theorem 1.1.} Under the assumptions of Theorem 1.1, by Theorem 4.1, we know that (1.1) is solvable in the
degenerate elliptic region $(t, x) \in [-T_0, 0] \times \R^n$, moreover, the corresponding solution $u(t,x)$  satisfies
$$u(t,x)\in C\Big([-T_0, 0], L^{p_0}(\Bbb R^n)   \Big), \ \ \p_t u(t,x)\in C\Big([-T_0, 0], L^{p_1}(\Bbb R^n)\Big).  \eqno{(6.1)}$$
On the other hand, by the initial data $(u(0,x), \p_tu(0,x))$ given in (6.1), it follows from Theorem 5.1 that the problem (1.1) admits a unique solution in the degenerate hyperbolic region $(t,x)
\in [0, T_0]\times\Bbb R^n$ for small constant $T_0>0$, which satisfies
$$u(t,x)\in C([0, T_0], L^{p_0}(\Bbb R^n)). \eqno{(6.2)} $$
Therefore, combining (6.1) with (6.2) yields Theorem 1.1.\hfill $\square$

\vskip 0.4 true cm

\end{document}